\documentclass[11pt]{article}

\usepackage[english]{babel} 
\usepackage[T1]{fontenc} 
\usepackage{amsmath,amsfonts,amssymb,amsthm} 

\usepackage{cases} 
 \usepackage{bm} 
\usepackage{enumitem}

\usepackage{xcolor} 
\usepackage{graphicx} 
\usepackage{pdfpages} 
\usepackage{subfig}
\usepackage{overpic}
\usepackage{tikz}
\usetikzlibrary{decorations.pathreplacing}
\tikzset{%
  show curve controls/.style={
    postaction={
      decoration={
        show path construction,
        curveto code={
          \draw [blue] 
            (\tikzinputsegmentfirst) -- (\tikzinputsegmentsupporta)
            (\tikzinputsegmentlast) -- (\tikzinputsegmentsupportb);
          \fill [red, opacity=0.5] 
            (\tikzinputsegmentsupporta) circle [radius=.5ex]
            (\tikzinputsegmentsupportb) circle [radius=.5ex];
        }
      },
      decorate
}}}
\usetikzlibrary{shapes.misc}
\tikzset{cross/.style={cross out, draw=black, minimum size=2*(#1-\pgflinewidth), inner sep=0pt, outer sep=0pt},
cross/.default={2.5pt}}
\usepackage{booktabs} 
\usepackage{colortbl}

\usepackage{authblk}
\usepackage[nottoc]{tocbibind}
\usepackage[titletoc,title]{appendix}

\numberwithin{equation}{section}
\providecommand{\keywords}[1]{\textbf{\textit{Keywords:}} #1}
\definecolor{blun}{cmyk}{0.8, 0.5, 0, 0.7}

\definecolor{Lpos}{rgb}{0.91, 0.33, 0.5}
\definecolor{Lneg}{rgb}{0, 0.81, 0.82}

\newcommand{\Lap}{\Delta}
\newcommand{\tr}{\textnormal{tr}}
\newcommand{\N}{\mathbb N}

\DeclareMathOperator{\Ker}{Ker}

\theoremstyle{plain}

\let\oldbibliography\thebibliography
\renewcommand{\thebibliography}[1]{%
  \oldbibliography{#1}%
  \setlength{\itemsep}{-1.2mm}%
}

\usepackage[pagewise, mathlines]{lineno}  
\usepackage[bookmarks,colorlinks,pdfauthor={Cinzia Soresina},linkcolor=blun,citecolor=blun]{hyperref}

\begin{document}

\title{\Large{Hopf bifurcations in the full SKT model and where to find them}}
\author[1]{C. Soresina}
\affil[1]{\footnotesize{Institut f\"ur Mathematik und wissenschaftliches Rechnen, \break Karl--Franzens Universit\"at Graz, Heinrichstr. 36, 8010 Graz, Austria}}
\maketitle

\begin{abstract}
\noindent
In this paper, we consider the Shigesada--Kawasaki--Teramoto (SKT) model, which presents cross-diffusion terms describing competition pressure effects. Even though the reaction part does not present the activator--inhibitor structure, cross-diffusion can destabilise the homogeneous equilibrium. However, in the full cross-diffusion system and weak competition regime, the cross-diffusion terms have an opposite effect and the bifurcation structure of the system modifies increasing the interspecific competition pressure. The major changes in the bifurcation structure, the type of pitchfork bifurcations on the homogeneous branch, as well as the presence of Hopf bifurcation points are here investigated. Through weakly nonlinear analysis, we can predict the type of pitchfork bifurcation. Increasing the additional cross-diffusion coefficients, the first two pitchfork bifurcation points from super-critical become sub-critical, leading to the appearance of a multi-stability region. The interspecific competition pressure also influences the possible appearance of stable time-period spatial patterns appearing through a Hopf bifurcation point.
\end{abstract}

\keywords{bifurcations; cross-diffusion; SKT model; \texttt{pde2path}; Stuart--Landau; Hopf} 


\section{Introduction}\label{sec_int}

Spatial segregation describes a situation in which two competing species coexist, but they mainly concentrate in different regions of the habitat. It can be due to territoriality or aggressiveness of individuals, or it can emerge from competition for the same (limited) resource. Other possible outcomes of interspecific competition are competitive exclusion and niche differentiation. Territorial segregation leads to exclusive exploitation of the resources and it can minimize the encounters, and consequently also the conflicts, between individuals~\cite{wilson1975}. It has been observed in several animal species, such as birds, mammals, amphibians, fishes and insects, resulting in checkerboard patterns~\cite{diamond1975assembly}.

From a mathematical point of view, models able to capture this effect should describe a situation in which the abundance of a population reduces the presence of the competing one. Several mathematical models have been proposed to explain this effect. To describe the spatial distribution of the populations and their interactions, reaction--diffusion systems can be formulated, and spatial segregation can result from diffusion-driven instability and pattern-formation, through non-homogeneous solutions of the reaction--diffusion model. A two-species Lotka--Volterra--Gause reaction--diffusion model with constant diffusion coefficients constitutes the first attempt in this context, but it fails to produce such patterns, at least on convex domains~\cite{kishimoto1985spatial}. This suggests that spatial segregation in a two-species competition-diffusion system must be caused by other mechanisms, such as competition-pressure \cite{shigesada1979spatial, mimura1981stationary} or non-convexity or non-homogeneity of the habitat \cite{matano1983pattern, ei1990pattern, shigesada1979spatial, levin1974dispersion}.

To account for stable inhomogeneous steady states exhibiting spatial segregation, the so-called SKT model was proposed in~\cite{shigesada1979spatial}. In addition to standard diffusion terms, the model includes nonlinear self-diffusion and cross-diffusion terms, modelling additional movements of individuals due to intra- and interspecific competition pressure. Interspecific competition pressure (cross-diffusion) translated into nonlinear diffusion terms depending on the presence of the competing species. This particular form might seem artificial at a first sight. However, a cross-diffusion term appears in the fast-reaction limit of a three-species system presenting only standard diffusion, competition and fast-reaction terms modelling the switch between two states 
\cite{iida2006diffusion}. The case in which self-diffusion is neglected and only one cross-diffusion term is considered (namely only one population has the special ability to avoid the other one) is often called \emph{triangular case} or \emph{triangular cross-diffusion system}, while \emph{full cross-diffusion system} refers to the case in which both equations present cross-diffusion terms.

The SKT model has attracted a lot of attention from different points of view. The first main theoretical result on global existence and regularity of the time depending solutions was obtained in~\cite{Ama88, Ama90}, where a general theory about the existence of local solutions for general quasilinear parabolic PDEs can be found. In~\cite{GalGarJun03} an entropy structure was discovered, and then entropy-based methods were then generalized~\cite{Jun16,DevLepMouTre15}. The convergence of solutions of the fast-reaction system to the ones of the cross-diffusion system has been rigorously proven for the triangular case in \cite{desvillettes2015new}. The convergence of the bifurcation structure of the three-species and the four-species systems to the cross-diffusion one has been investigated in \cite{iida2006diffusion, CKCS} and in \cite{CKetal6} respectively. For more details, we refer to \cite{breden2021influence, CKCS} and the references therein.

Starting from the seminal paper~\cite{shigesada1979spatial}, the question of the existence of non-homogeneous steady states for the SKT system, when cross-diffusion terms are taken into account, has been extensively investigated, both numerically and theoretically. In this regard, we must distinguish two different regimes: weak and strong competition. In the weak competition regime, the homogeneous system (when all diffusion terms are neglected) admits unstable non-coexistence equilibria and a stable coexistence one. With only standard diffusion, in a convex domain and with zero-flux boundary conditions, any non-negative solution generically converges to the coexistence steady state, and this implies that the two species coexist but their densities are homogeneous in the whole domain~\cite{kishimoto1985spatial}. In the strong competition case instead, for the homogeneous system, the coexistence steady state is unstable, while the non-coexistence ones are stable. Adding only standard diffusion, in a convex domain and with zero-flux boundary conditions, it has been shown that if positive and non-constant steady states exist, they must be unstable~\cite{kishimoto1985spatial}, and numerical simulations suggest that any non-negative solution generically converges to one of the two (non-trivial) non-coexistence states, predicting the competitive exclusion of the two species.

In the weak competition regime, the surprising effect is that if one of the cross-diffusion terms is sufficiently large compared to all other parameters, then the homogeneous coexistence steady state loses its stability and non-homogeneous steady states appear~\cite{LouNi96, LouNiYot04, ni2014existence, lou2015pattern, mori2018numerical}. This is known as cross-diffusion driven instability, being the cross-diffusion term the key ingredient that destabilises the homogeneous equilibrium. Besides, the shape and the amplitude of these patterns can be predicted~\cite{gambino2012turing}. We refer to ~\cite{IidNimYam18, Jun10, breden2021influence} and the references therein for a broader discussion about both regimes.

Of particular interest are the changes of steady states, and more generally of the bifurcation structure, under parameter variation. In addition to theoretical methods and criteria, the bifurcation structure of steady states discloses the behaviour of solutions far from homogeneous branch. For classical reaction--diffusion systems on bounded domains (and in dimension~$1,2,3$), the bifurcation structure can be numerically computed with, for instance, \texttt{pde2path}~\cite{uecker2021numerical, uecker2014pde2path, uecker2021continuation}. This package is an advanced continuation/bifurcation software based on the FEM discretization of the stationary elliptic problem exploiting the package \texttt{OOPDE}~\cite{prufert2014oopde} for the FEM discretization. Since the software is quite flexible, it has also been used beyond its standard-setting, for instance, to treat cross-diffusion systems~\cite{CKCS}, and spectral fractional diffusion \cite{ehstand2021numerical}. Moreover, thanks to powerful \emph{computer-assisted techniques} \cite{breden2013global,breden2018existence} developed in the last three decades, and recently extended to treat non-linear diffusion terms \cite{breden2021computer}, the approximated solutions found with \texttt{pde2path} can be validated rigorously a posteriori.

Thanks to the interplay between linearised analysis and numerical continuation, in \cite{breden2021influence} the full SKT model has been investigated over a large range of parameters, revealing some interesting effects of the cross-diffusion terms on the steady states. In particular, from the linearised analysis, we have that the cross-diffusion terms have an opposite effect in the destabilisation of the homogeneous equilibrium. However, this holds close to the homogeneous branch, while the investigation of how the bifurcation structure modifies with respect to the cross-diffusion coefficients revealed the presence of multistability regions, Hopf bifurcation points and also the presence of stable non-homogeneous solutions beyond the usual parameter range of investigation, suggesting that their influence is more complex and rich than the one predicted by the linearised analysis.

Finally, in addition to steady states, it has been proven that the SKT model can also exhibit other types of patterns, such as stable time-periodic solutions arising through a Hopf bifurcation point. Such solutions describe a dynamic coexistence between the two species. In the strong competition regime, the existence of periodic solutions has been proven in~\cite{kan1993stability}. The existence of stable time-periodic solutions which bifurcate from a Hopf bifurcation point in the weak competition regime has been proven in~\cite{izuhara2018spatio}, in the triangular case applying the center manifold theory and the standard normal form theory. The important aspect highlighted here is the presence of a doubly degenerate point at the intersection of the neutral stability curves related to the 1- and 2-modes.

This work aims to further analyse the full cross-diffusion SKT model, namely with a cross-diffusion term in both equations, in the weak competition regime. We are interested in studying the influence of the additional interspecific-competition pressure (cross-diffusion) on the bifurcation structure, extending the study carried out in \cite{breden2021influence}. In particular, we focused on the type of pitchfork bifurcation on the homogenous branch (related to the stability of the bifurcating branches close to the homogeneous one) and on the presence of Hopf bifurcation points on the bifurcating branch corresponding to the 1-mode.

It has been pointed out in \cite{breden2021influence} that the model presents multistability of solutions, meaning that a particular stable inhomogeneous solution can coexist with the homogeneous one for suitable parameter values. This is an important aspect in ecology, since in this case by perturbing the system it is possible to pass from homogeneous distribution of the species on the habitat to spatial segregation. From the mathematical point of view, this situation is related to the type of pitchfork bifurcations (sub- or super-critical) on the homogenous branch, and its dependence on the cross-diffusion coefficients. Therefore, we aim to  analytically characterise these bifurcation points. This effect cannot be captured using linearised analysis only, but it can be achieved through weakly nonlinear analysis, deriving the Stuart--Landau equation at the bifurcation point, exploiting the technique presented in \cite{gambino2012turing, gambino2013pattern, gambino2016super}. Once the coefficients that characterise the pitchfork bifurcation are obtained, it is possible to study their dependence on the cross-diffusion coefficients, varying the cross-diffusion parameters. Thanks to this result, we can better characterize the effect of the cross-diffusion term on the bifurcation structure and determine the parameter regions in which multistability appears.

Another important aspect of the bifurcation structure of the SKT system is the possible presence of Hopf bifurcation points, suggesting the formation of time-periodic spatial patterns. The influence of the cross-diffusion terms on the Hopf points, as well as on the effective presence, type and stability properties of these time-varying patterns, are biologically relevant. To investigate the possible scenarios, the analytical results obtained in \cite{izuhara2018spatio} can be extended to the full cross-diffusion case, and combined with the weakly nonlinear analysis and numerical continuation results.

Thanks to its interplay between linearised analysis, weakly nonlinear analysis and numerical continuation, this work constitutes a step forward in the analytical understanding of the bifurcation structure of the SKT system.

The paper is organized as follow. In Section \ref{sec_model} the model is introduced and the linearised analysis needed in the following is reported. In Section \ref{sec_bif} the major changes in the bifurcation structure are presented. Section \ref{sec_SL} is devoted to the weakly nonlinear analysis and the study of the type of pitchfork bifurcation on the homogeneous branch, while in Section \ref{sec_Hopf} the presence of Hopf bifurcation points is investigated. Finally, in Section \ref{sec_concl} some concluding remarks can be found. The Matlab files needed by \texttt{pde2path}, as well as the Matlab scripts related to Sections \ref{sec_SL} and \ref{sec_Hopf} can be found on GitHub \cite{GitHubSKT, GitHubSL}. They can be used to check and reproduce the numerical results, as a tutorial or as a starting point for further investigations.

\section{The full SKT model and linear stability analysis}\label{sec_model}
We consider here the so-called SKT model,  proposed in~\cite{shigesada1979spatial} to account for stable inhomogeneous steady states exhibiting spatial segregation of two competing species.  We denote with $u(t,x),\; v(t,x)\geq 0$ the population densities of two species at time $t$ and position $x$, confined and competing for resources on a bounded and connected domain $\Omega~\subset~\mathbb{R}^N$.  The system describing the dynamics writes
\begin{equation}
\begin{cases}
\partial_t u=\Lap((d_1+d_{11} u +d_{12} v)u)+(r_1-a_1 u-b_1 v)u,& \textnormal{on } \mathbb{R}_+\times \Omega,\\[0.2cm]
\partial_t v=\Lap((d_2+d_{22} v +d_{21} u)v)+(r_2-b_2 u-a_2 v)v,& \textnormal{on }  \mathbb{R}_+\times \Omega,
\vspace{0.2cm}
\\
\dfrac{\partial}{\partial n}u=\dfrac{\partial}{\partial n}v=0,&\textnormal{on } \mathbb{R}_+\times \partial\Omega,\\[0.2cm]
u(0,x)=u_{in}(x),\; v(0,x)=v_{in}(x),& \textnormal{on } \Omega,
\end{cases}\label{cross}
\end{equation}
where  the coefficients $r_i,\,a_i,\,b_i \, (i=1,2)$ are the intrinsic growth, the intraspecific competition and the interspecific competition rates.  Parameters $d_1$ and $d_2$ describes the diffusion, while $d_{11},\,d_{22}$ and $d_{12},\,d_{21}$ stand for competition pressure,  and are called self- and cross-diffusion coefficients.  To avoid confusion,  we will refer to $d_1$ and $d_2$ as the \emph{standard} diffusion coefficients.  Throughout this paper we consider the cross-diffusion system \eqref{cross} and assume that the standard diffusion coefficients are positive, and that all the other coefficients are non-negative.  

Since we are mainly interested in the influence of cross-diffusion coefficient in the appearance of non-homogeneous steady states, we consider $d_1=d_2=:d$, as already done in previous studies. Moreover,  it is known that self-diffusion coefficients inhibit the emergence of those type of stationary solution (see \cite[Section 5]{breden2021influence}),  so we consider here $d_{11}=d_{22}=0$.

Looking at homogeneous steady states, system \eqref{cross} admits the total extinction $(0,0)$, two non-coexistence states $(\bar u,0)~=~(r_1/a_1,0)$ and $(0,\bar v)=(0,r_2/a_2)$, and one coexistence state 
\begin{equation*}(u_*,v_*)=\left( \dfrac{r_1a_2-r_2b_1}{a_1a_2-b_1b_2},\dfrac{r_2a_1-r_1b_2}{a_1a_2-b_1b_2}\right).
\end{equation*}
The non-coexistence equilibria exist for all the parameter values,  while the admissibility of the coexistence steady state (i.e.~positivity) leads to two regimes, the weak competition and the strong competition regimes \cite{breden2021influence}.  In this paper we are interested in the weak competition case,  namely when 
\begin{equation}\label{weak_comp}
\dfrac{b_1}{a_2}<\dfrac{r_1}{r_2}<\dfrac{a_1}{b_2}.
\end{equation}
In this case, for the homogeneous system (when all diffusion coefficients are taken equal to zero), the coexistence steady state is stable, while the non-coexistence ones are unstable. With only standard diffusion, in a convex domain and with zero-flux boundary conditions, any non-negative solution generically converges to the coexistence steady state $(u_*,v_*)$, and this implies that the two species coexist but their densities are homogeneous in the whole domain~\cite{kishimoto1985spatial}. 

We start by studying ``mode by mode'' the linear stability of $(u_*,v_*)$~\cite{Henry,KuehnBook1},  and we consider the eigenfunctions $\psi_k$ and associated eigenvalues $-\lambda_k$ of the Laplacian with zero Neumann boundary conditions, which satisfy $\lambda_0=0$, $\lambda_k>0$ for all $k\in\N,\,{k\geq 1}$, and $\lambda_k\rightarrow+\infty$ as $k\to+\infty$ (we always assume that the eigenvalues are labeled in ascending order). 

The Jacobian matrix of the reaction part and the linearisation of the diffusion part of~\eqref{cross}, evaluated at the equilibrium $(u_*,v_*)$, are
\begin{equation*}
K= \begin{pmatrix}-a_1u_* & -b_1u_*\\-b_2v_*& -a_2v_*\end{pmatrix}, \quad
D=\begin{pmatrix}d+d_{12}v_*& d_{12}u_*\\d_{21}v_* & d+d_{21}u_*\end{pmatrix}.
\end{equation*}
Note that the coexistence equilibrium is stable in the weak competition regime, hence we have $\tr K<0$ and in particular $\det K>0$ (the inter- and intraspecific competition rates~$\;a_i,\, b_i,\; i=1,2$ are fixed parameters).
Then, the characteristic matrix associated to the $k$-th mode, $k\in\N$, is
\begin{equation}\label{matrixM_k}
M_k=K-\lambda_kD=
\begin{pmatrix}
-a_1u_*-(d+d_{12}v_*)\lambda_k   &  -b_1u_*-d_{12}u_*\lambda_k\\
-b_2v_* -d_{21}v_*\lambda_k         & -a_2v_*-(d+d_{21}u_*)\lambda_k
\end{pmatrix},
\end{equation}
and its determinant can be written as a second order polynomial in the bifurcation parameter $d$
\begin{equation}\label{detM*d}
P_k(d):=\det M_k=\lambda_k^2 d^2 +
(d_{12}v_*\lambda_k^2+d_{21}u_*\lambda_k^2-\tr K\lambda_k)d
-d_{12}\alpha \lambda_k- d_{21}\beta \lambda_k +\det K,
\end{equation}
where 
\begin{equation}\label{eq:alphabeta}
\alpha:=(b_2u_*-a_2v_*)v_*, \qquad 
\beta:=(b_1v_*-a_1u_*)u_*, \qquad 
\det K=(a_1a_2-b_1b_2)u_*v_*.
\end{equation}
Since $\det K$ is always positive and the trace of $M_k$ is always negative, therefore the $k$-th mode is stable if $P_k(d)>0$, unstable if $P_k(d)<0$, and a bifurcation occurs at the critical value $d_c$ of the bifurcation parameter for $P_k(d_c)=0$.  Introducing
\begin{equation}
\label{ABC}
A_k=\lambda_k^2,\quad B_k=d_{12}v_*\lambda_k^2+d_{21}u_*\lambda_k^2-\tr K\lambda_k, \quad C_k=-d_{12}\alpha \lambda_k- d_{21}\beta \lambda_k +\det K,
\end{equation}
the critical value $d_c$ can be easily computed as the positive solution to 
\begin{equation*}
P_k(d)=A_k d^2 + B_k d + C_k=0.
\end{equation*}
Note that the critical value $d_c$ depends on the eigenvalue $\lambda_k$ considered, and that no bifurcation can happen for $k=0$, since $P_0$ reduces to $\det K$, which is independent of $d$.
Obviously $A_k\geq 0$ and $B_k>0$. Therefore, a bifurcation associated to the $k$-th mode ($k\geq 1$) can occur if and only if $C_k<0$.  In this case, we can compute the critical value $d_c=d_c(\lambda_k, d_{12},d_{21})$ as
\begin{equation}\label{formula_d_c}
d_c=\dfrac{-B_k+\sqrt{B_k-4A_kC_k}}{4A_k},
\end{equation}
being the discriminat always positive.  The signs of $\alpha$ and $\beta$, which depend on the parameter values $r_i,\,a_i,\, b_i,\, (i=1,2),$ are thus crucial, as they change the monotonicity of $C_k$ with respect to $d_{12}$ and $d_{21}$ respectively.  The study of the possible combination with respect to the parameter values $r_i,\,a_i,\, b_i,\, (i=1,2)$ highlights the opposite role of the cross-diffusion cofficients $d_{12}$ and $d_{21}$ in destabilising the homogeneous equilibrium $(u_*,v_*)$ \cite{breden2021influence}.  

\begin{figure}[!ht]
\centering
\subfloat[$d_{21}=0$ \label{nsc_d21_0}]{
\begin{overpic}[width=0.45\textwidth,tics=10]{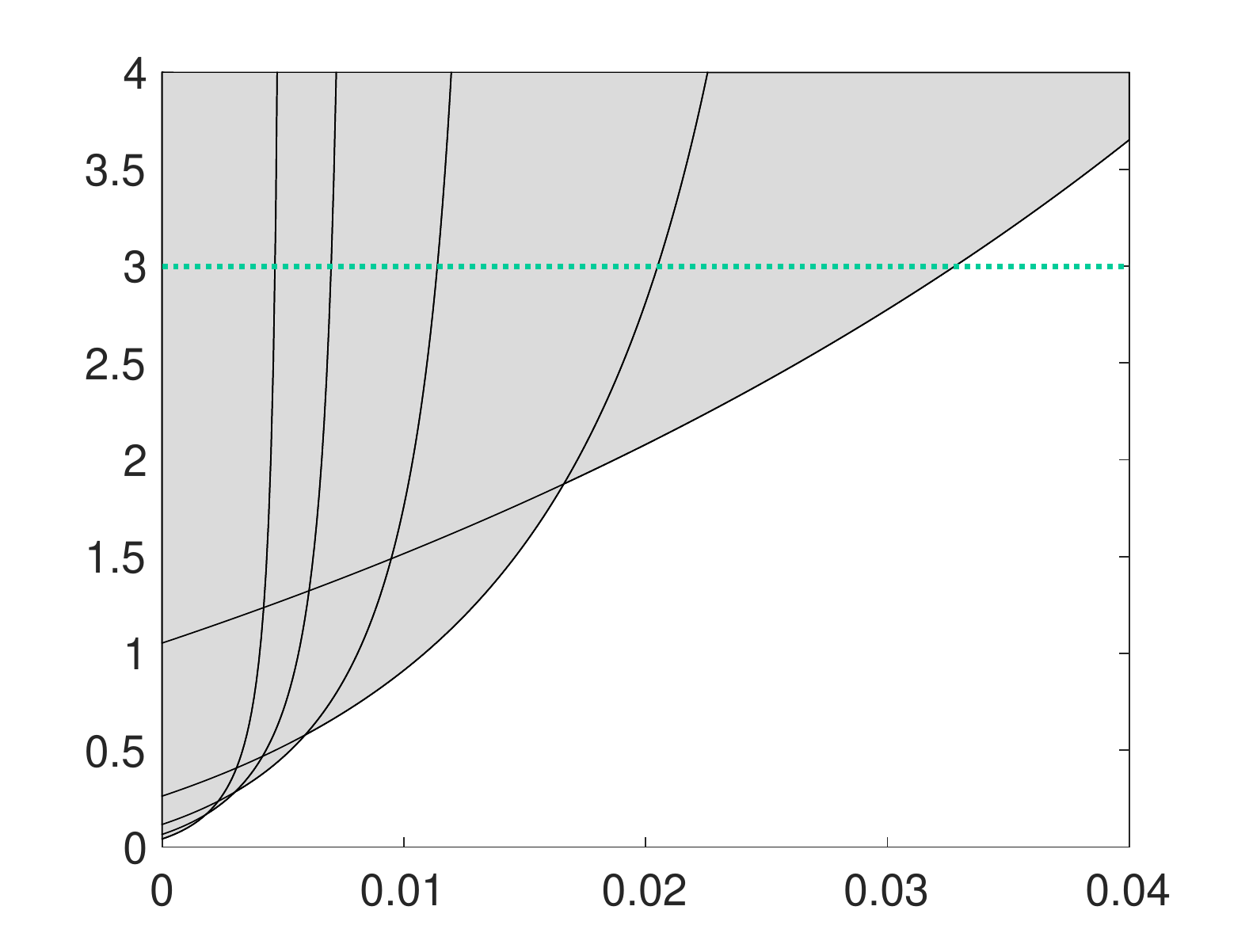}
\put(-2,60){$d_{12}$}
\put(90,-2){$d$}
\put(45,30){$(\hat{d},\hat{d}_{12})$}
\end{overpic} 
}\hspace{0.1cm}
\subfloat[$d_{21}=0.025$\label{nsc_d21_0p025}]{
\begin{overpic}[width=0.45\textwidth,tics=10]{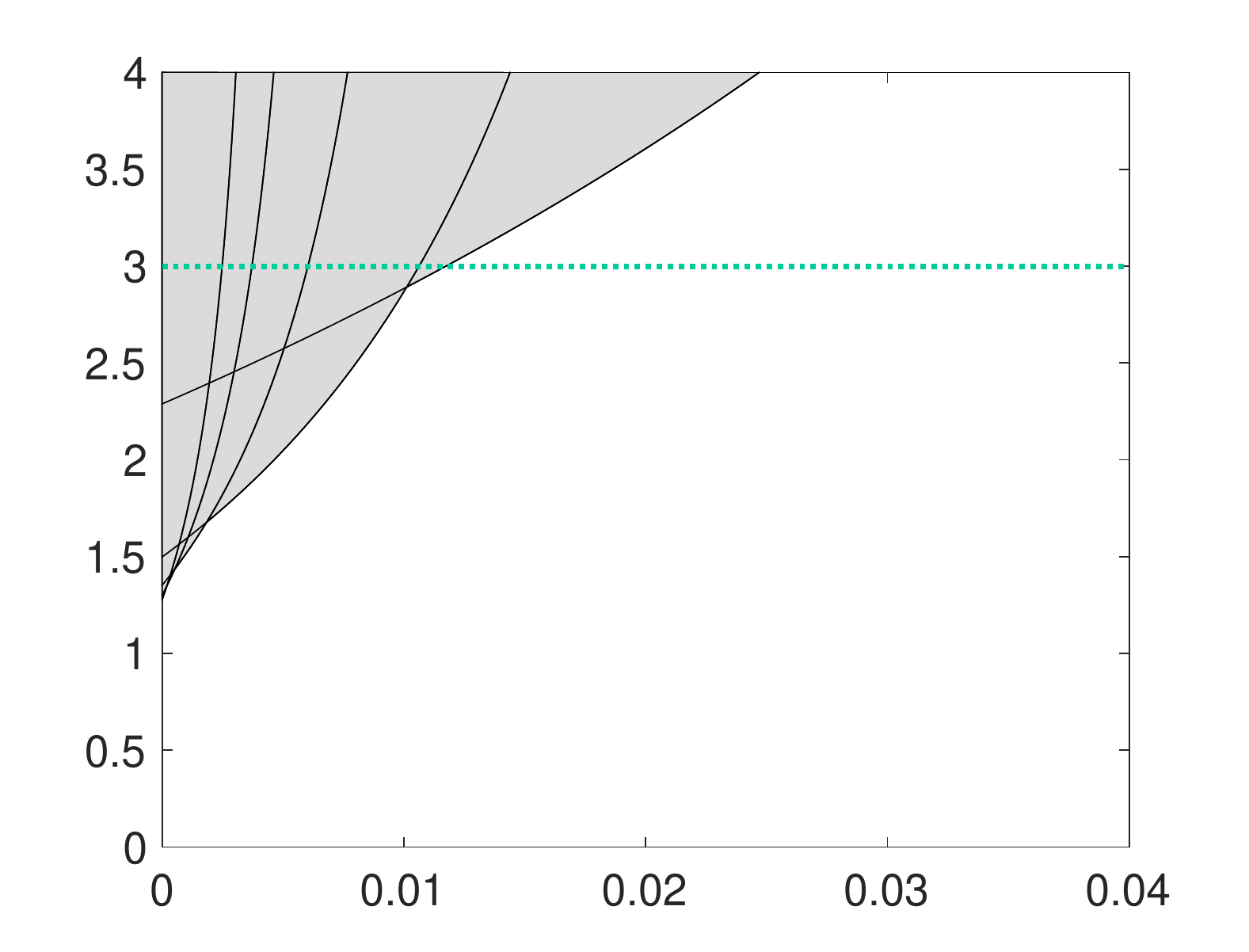}
\put(-2,60){$d_{12}$}
\put(90,-2){$d$}
\put(30,45){$(\hat{d},\hat{d}_{12})$}
\end{overpic} 
}
\caption{Neutral stability curves for $\lambda_k,\,k=1,\dots,5$ for different values of the cross-diffusion coefficient $d_{21}$. The white area denotes the stability of the homogenous steady-state $(u_*,v_*)$, while in the grey region the homogeneous steady state is destabilized and stable non-homogeneous stationary solutions appear. The green dotted horizontal line marks the ``usual'' value of parameter $d_{12}$.  The doubly-degenerate point $(\hat{d},\hat{d}_{12})$ corresponds to the intersection of the neutral stability curves associated to the 1- and 2-modes.}\label{neutral_stab_curves}
\end{figure}
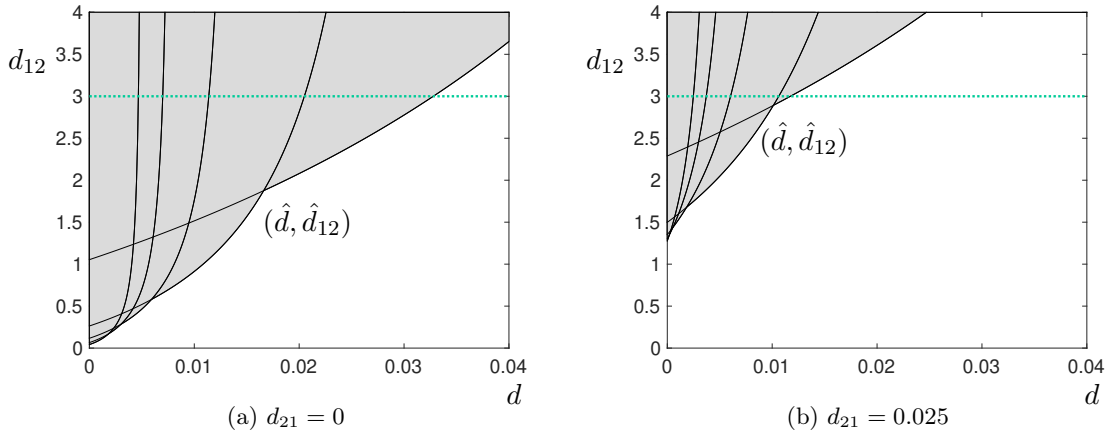

This can also been observed by looking at the neutral stability curves
\begin{equation}
\label{nsc}
d_{12}(d,d_{21},\lambda_k)=\dfrac{\lambda_k^2d^2+d_{21}du_*\lambda_k^2-d\lambda_k\tr K-d_{21}\beta\lambda_k+\det K}{\alpha\lambda_k-dv_*\lambda_k^2}, \qquad k\geq 1.
\end{equation}
In Figure \ref{neutral_stab_curves}, these curve for different $\lambda_k,\,k=1,\dots,5$ are shown, considering the ``usual'' setting~\cite{iida2006diffusion, izuhara2008reaction, breden2018existence, CKCS, breden2021influence} reported in Table \ref{tab:param} for the reader's convenience.

\begin{table}
\centering
\begin{tabular}{ccccccc>{\columncolor[gray]{0.9}}c>{\columncolor[gray]{0.9}}c}
$\Omega$&$r_1$&$r_2$&$a_1$&$a_2$&$b_1$&$b_2$&$\alpha$&$\beta$\\
\midrule[2pt]
(0,1)&5&2&3&3&1&1&$+$&$-$\\
\bottomrule[2pt]
\end{tabular}
\caption{The parameter sets used in the numerical simulations, relevant to the weak competition regime. The sign of the quantity $\alpha$ and $\beta$ in \eqref{eq:alphabeta} is reported. }
\label{tab:param}
\end{table}

Figure \ref{nsc_d21_0} belongs to the triangular case $d_{21}=0$, while in Figure \ref{nsc_d21_0p025} the value of $d_{21}$ has been increased. The white area corresponds to the stable region of the homogeneous steady state $(u_*,v_*)$, while in the grey region the homogeneous steady state is destabilized and stable non-homogeneous stationary solutions appear. We can observe that the grey region reduces (when $d_{12}$ is fixed) as the cross-diffusion coefficients $d_{21}$ increases, since the bifurcation points (points of the neutral stability curves for a fixed value $d_{12}$) move towards zero.

We can also identify the presence of doubly degenerate point $(\hat{d},\hat{d}_{12})$ where two curves intersect, meaning that the corresponding two bifurcation points on the homogeneous branch switch their position. It has been proven for the triangular system that in the vicinity of this point a Hopf bifurcation appears for $d_{12}<\hat{d}_{12}$ on the bifurcation branch corresponding to $\lambda_1$ \cite{izuhara2018spatio}, and the system shows stable time-periodic non-homogeneous solutions. In the full cross-diffusion system, the cross-diffusion coefficient $d_{21}$ shifts up the curves. Then,  in the full cross-diffusion system ($d_{21}>0$) the doubly degenerate point persists and, for the usual parameter set, happens for greater values of $d_{12}$ than for the triangular case.  This suggests
the possible presence of a Hopf point in the bifurcation structure of the full-cross diffusion system. 

\section{Bifurcation diagrams}\label{sec_bif}
As observed in the previous section, the effect of the cross-diffusion coefficient $d_{21}$ on the bifurcation point of the homogeneous branch can be studied by looking at equations \eqref{detM*d}. In particular, in the weak competition case, it is known that the cross-diffusion coefficients have an opposite effect in the destabilisation of the homogeneous equilibrium \cite{gambino2012turing}, and it has been studied in detail in \cite{breden2021influence}. However, how the bifurcation structure (in particular far from the homogeneous branch) behaves and modifies for increasing values of $d_{21}$ is not predictable through linearised analysis only. Nethertheless, it can be numerically computed and to this end, we exploit the continuation software for PDEs \texttt{pde2path}~\cite{dohnal2014pde2path, uecker2018hopf, uecker2014pde2path}, based on a FEM discretization of the stationary problem, and the suitable setting for cross-diffusion terms proposed in \cite{CKCS}. We refer to \cite{GitHubSKT} for the Matlab scripts needed to compute the bifurcation structure in \texttt{pde2path}.

As in the previous section, we consider the usual setting, such as the 1D domain $\Omega=(0,1)$ and the parameter set already used in \cite{iida2006diffusion, izuhara2008reaction, breden2018existence, CKCS,breden2021influence}, reported in Table \ref{tab:param}. We fix the cross-diffusion coefficient $d_{12}=3$. In Figure \ref{BD_d21} we show the bifurcation diagrams obtained for increasing values of the cross-diffusion parameter $d_{21}$.  In particular, the first diagram corresponds to the triangular case ($d_{21}=0$), while the others belong to the full cross-diffusion system. The bifurcation parameter is the standard diffusion coefficient $d$, while in the $y$-axis the quantity $v(0)$ is reported. Thicker lines in the bifurcation diagrams denote stable solutions, while we use thinner lines for unstable ones. Circles, crosses and diamonds indicate the presence of branch points, fold points and Hopf points, respectively. The black horizontal line denotes the homogeneous branch (the homogeneous steady state $(u_*,v_*)$ is independent of $d$), while the coloured branches bifurcating from the homogeneous one correspond to different eigenvalues $\lambda_k$. In particular, the blue, red, green, yellow branches correspond to $\lambda_i,\, i=1,2,3,4$, respectively, while the magenta branch originates from a secondary bifurcation point. We also observe that the blue and the red branches ``intersects'' at a secondary bifurcation point, at which the red lower branch becomes stable.

As predicted from the linearised analysis, the bifurcation points on the homogeneous branch move to the left as $d_{21}$ increases. Consequently, the whole bifurcation structure tends to shrink to the left. We observe that, while the bifurcation structure in Figures \ref{BD_d21_0} and \ref{BD_d21_0p01} are qualitatively similar, from Figure \ref{BD_d21_0p02} the first bifurcating branch (in blue) undergoes to a major deformation: the first bifurcation point from super-critical becomes sub-critical. This leads to a multi-stability region in which the system presents both homogeneous and non-homogeneous stable solutions. Increasing the value of $d_{21}$ even further, other qualitative changes can be observed. We observe a further change in the type of pitchfork bifurcation (super-/sub-critical) of the second bifurcation point, and the change in stability property of the red branch follows thanks to a fold bifurcation. Moreover, on the first (blue) bifurcation branch the continuation software detects a Hopf bifurcation point, related to a pair of imaginary eigenvalues. The Hopf bifurcation point persists for greater values of the bifurcation parameters. 

\begin{figure}
\centering
\subfloat[\label{noH}]{
\begin{tikzpicture}
\draw (-0.5,0) -- (3,0);
\draw[ultra thick] (3,0) -- (3.5,0);
\draw[ultra thick] (3,0) ..  controls ++(100:1) and ++(90:1)   .. (0.2,0.8);
\draw (0.2,0.8) ..  controls ++(-90:0.5) and ++(90:0.7) .. (0.5, -0.4)
..  controls ++(-90:0.2) and ++(90:0.2) .. (0.2, -0.8);
\draw[ultra thick] (3,0) ..  controls ++(-100:1) and ++(-90:1)   .. (0.2,-0.8);
\draw (1,0) ..  controls ++(100:0.3) and ++(-30:1)   .. (0,1.5);
\draw[ultra thick] (0,1.5) ..  controls ++(150:0.2) and ++(-20:0.2) .. (-0.5, 1.75) ;
\draw (1,0) ..  controls ++(-100:0.3) and ++(30:0.1)   .. (0.5,-0.4);
\draw[ultra thick] (0.5,-0.4) ..  controls ++(-160:0.2) and ++(10:0.2) .. (-0.5, -0.7);
\draw[fill=white] (1,0) circle (0.9mm);
\draw[fill=white] (3,0) circle (0.9mm);
\draw[fill=white] (0,1.5) circle (0.9mm);
\draw[fill=white] (0.5,-0.4) circle (0.9mm);
\draw (.2,.8) node[cross] {};
\draw (.2,-.8) node[cross] {};
\draw (3.3,0.15) node[] {{\color{gray}\tiny$0$}};
\draw (2.8,.7) node[] {{\color{gray}\tiny$0$}};
\draw (2.8,-.7) node[] {{\color{gray}\tiny$0$}};
\draw (2.7,.15) node[] {{\color{gray}\tiny$1$}};
\draw (0.4,.4) node[] {{\color{gray}\tiny$1$}};
\draw (0.5,-.7) node[] {{\color{gray}\tiny$1$}};
\draw (1,.4) node[] {{\color{gray}\tiny$1$}};
\draw (1,-.3) node[] {{\color{gray}\tiny$1$}};
\draw (0,.15) node[] {{\color{gray}\tiny$2$}};
\draw (-0.2,-0.4) node[] {{\color{gray}\tiny$0$}};
\draw (-.2,1.8) node[] {{\color{gray}\tiny$0$}};
\end{tikzpicture}
}
\subfloat[\label{preH}]{
\begin{tikzpicture}
\draw (-0.5,0) -- (2,0);
\draw[ultra thick] (2,0) -- (3.5,0);
\draw[ultra thick] (2.3,0.5) ..  controls ++(90:1) and ++(90:1)   .. (0.2,0.8);
\draw (0.2,0.8) ..  controls ++(-90:0.5) and ++(90:0.7) .. (0.5, -0.4)
..  controls ++(-90:0.2) and ++(90:0.2) .. (0.2, -0.8);
\draw (2.3,0.5) ..  controls ++(-90:0.2) and ++(90:0.2) .. (2, 0)
..  controls ++(-90:0.2) and ++(90:0.2) .. (2.3, -0.5);
\draw[ultra thick] (2.3,-0.5) ..  controls ++(-100:1) and ++(-90:1)   .. (0.2,-0.8);
\draw (1,0) ..  controls ++(100:0.3) and ++(-30:1)   .. (0,1.5);
\draw[ultra thick] (0,1.5) ..  controls ++(150:0.2) and ++(-20:0.2) .. (-0.5, 1.75) ;
\draw (1,0) ..  controls ++(-100:0.3) and ++(30:0.1)   .. (0.5,-0.4);
\draw[ultra thick] (0.5,-0.4) ..  controls ++(-160:0.2) and ++(10:0.2) .. (-0.5, -0.7);
\draw[fill=white] (0,1.5) circle (0.9mm);
\draw[fill=white] (0.5,-0.4) circle (0.9mm);
\draw[fill=white] (1,0) circle (0.9mm);
\draw[fill=white] (2,0) circle (0.9mm);
\draw (.2,.8) node[cross] {};
\draw (.2,-.8) node[cross] {};
\draw (2.3,.5) node[cross] {};
\draw (2.3,-.5) node[cross] {};
\draw (2.3,.15) node[] {{\color{gray}\tiny$0$}};
\draw (2.4,.8) node[] {{\color{gray}\tiny$0$}};
\draw (2.4,-.8) node[] {{\color{gray}\tiny$0$}};
\draw (2.05,.3) node[] {{\color{gray}\tiny$1$}};
\draw (2.05,-.3) node[] {{\color{gray}\tiny$1$}};
\draw (1.7,.15) node[] {{\color{gray}\tiny$1$}};
\draw (0.4,.4) node[] {{\color{gray}\tiny$1$}};
\draw (0.5,-.7) node[] {{\color{gray}\tiny$1$}};
\draw (1,.4) node[] {{\color{gray}\tiny$1$}};
\draw (1,-.3) node[] {{\color{gray}\tiny$1$}};
\draw (0,.15) node[] {{\color{gray}\tiny$2$}};
\draw (-0.2,-0.4) node[] {{\color{gray}\tiny$0$}};
\draw (-.2,1.8) node[] {{\color{gray}\tiny$0$}};
\end{tikzpicture}
}
\subfloat[\label{withH}]{
\begin{tikzpicture}
\draw (-0.5,0) -- (2,0);
\draw[ultra thick] (2,0) -- (3.5,0);
\draw (2.3,0.5) ..  controls ++(90:1) and ++(90:1)   .. (0.5,0.8);
\draw (0.5,0.8) ..  controls ++(-90:0.5) and ++(140:0.7) .. (1.15, -0.3)
..  controls ++(-90:0.2) and ++(90:0.2) .. (0.5, -0.8);
\draw (2.3,0.5) ..  controls ++(-90:0.2) and ++(90:0.2) .. (2, 0)
..  controls ++(-90:0.2) and ++(90:0.2) .. (2.3, -0.5);
\draw (2.3,-0.5) ..  controls ++(-100:1) and ++(-90:1)   .. (0.5,-0.8);
\draw[ultra thick] (2.3,-0.5) ..  controls ++(-100:0.9) and ++(-35:0.33)   .. (0.8,-1.34);
\draw[ultra thick] (2.3,0.5) ..  controls ++(85:0.85) and ++(40:0.43)   .. (0.8,1.34);
\draw (1,0) ..  controls ++(80:0.2) and ++(90:-0.1)   .. (1.3,0.5)
..  controls ++(90:0.1) and ++(-30:1)   .. (0,1.5);
\draw[ultra thick] (0,1.5) ..  controls ++(150:0.2) and ++(-20:0.2) .. (-0.5, 1.75) ;
\draw (1,0) ..  controls ++(-80:0.3) and ++(-90:-0.1)   .. (1.3,-0.5)
..  controls ++(-90:0.1) and ++(10:0.2) .. (-0.5, -0.7);
\draw[ultra thick](1.3,-0.5) ..  controls ++(-90:0.1) and ++(10:0.2) .. (-0.5, -0.7);
\draw[fill=white] (1,0) circle (0.9mm);
\draw[fill=white] (2,0) circle (0.9mm);
\draw[fill=white] (0,1.5) circle (0.9mm);
\draw[fill=white] (1.1,-0.3) circle (0.9mm);
\draw (.5,.8) node[cross] {};
\draw (.5,-.8) node[cross] {};
\draw (2.3,.5) node[cross] {};
\draw (2.3,-.5) node[cross] {};
\draw (1.3,.5) node[cross] {};
\draw (1.3,-.5) node[cross] {};
\draw[fill=white, rotate around={45:(0.8,1.27)}] (0.8,1.27) rectangle ++(4pt,4pt);
\draw[fill=white, rotate around={45:(0.8,-1.41)}] (0.8,-1.41) rectangle ++(4pt,4pt);
\draw (2.3,.15) node[] {{\color{gray}\tiny $0$}};
\draw (1.7,.15) node[] {{\color{gray}\tiny $1$}};
\draw (0,.15) node[] {{\color{gray}\tiny $2$}};
\draw (2.4,.8) node[] {{\color{gray}\tiny $0$}};
\draw (2.4,-.8) node[] {{\color{gray}\tiny $0$}};
\draw (2.05,.3) node[] {{\color{gray}\tiny $1$}};
\draw (2.05,-.3) node[] {{\color{gray}\tiny $1$}};

\draw (0.4,.4) node[] {{\color{gray}\tiny $1$}};
\draw (0.5,-1.2) node[] {{\color{gray}\tiny $2$}};

\draw (0.8,-0.4) node[] {{\color{gray}\tiny $1$}};
\draw (0.55,1.35) node[] {{\color{gray}\tiny $2$}};

\draw (1,.3) node[] {{\color{gray}\tiny $2$}};
\draw (1.2,-.13) node[] {{\color{gray}\tiny $2$}};
\draw (1.3,-.3) node[] {{\color{gray}\tiny $1$}};
\draw (1.3,0.75) node[] {{\color{gray}\tiny $1$}};

\draw (1.3,-0.7) node[] {{\color{gray}\tiny $0$}};
\draw (-.2,1.8) node[] {{\color{gray}\tiny $0$}};
\end{tikzpicture}
}
\caption{Qualitative representation of the bifurcation structure at the first and second bifurcation points.  Numbers along the branches indicate the number of eigenvalues with positive real part detected by the continuation software \texttt{pde2path}.}
\label{fig:BDschizzo}
\end{figure}
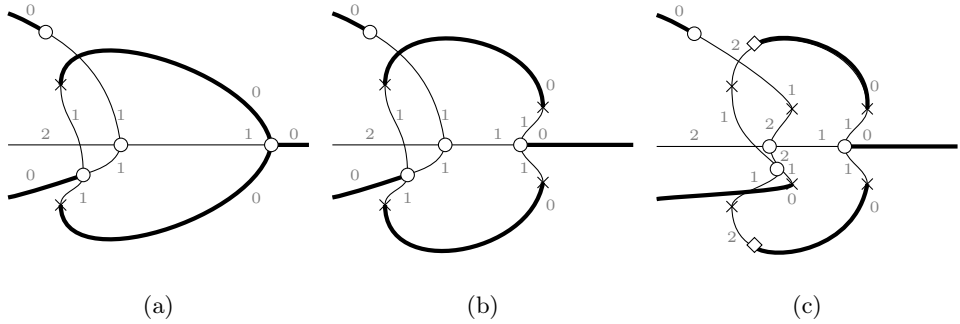

A qualitative representation of the first two branches is shown in Figure \ref{fig:BDschizzo}, where we also indicate the number of eigenvalues with positive real part detected numerically. Note that in Figure \ref{withH}, the Hopf bifurcation is not necessary to match the number of instabilities along the branch, but also a change of stability due to a fold point is compatible.

\begin{figure}
\centering
\subfloat[$d_{21}=0$\label{BD_d21_0}]{
\begin{overpic}[width=0.45\textwidth,tics=10]{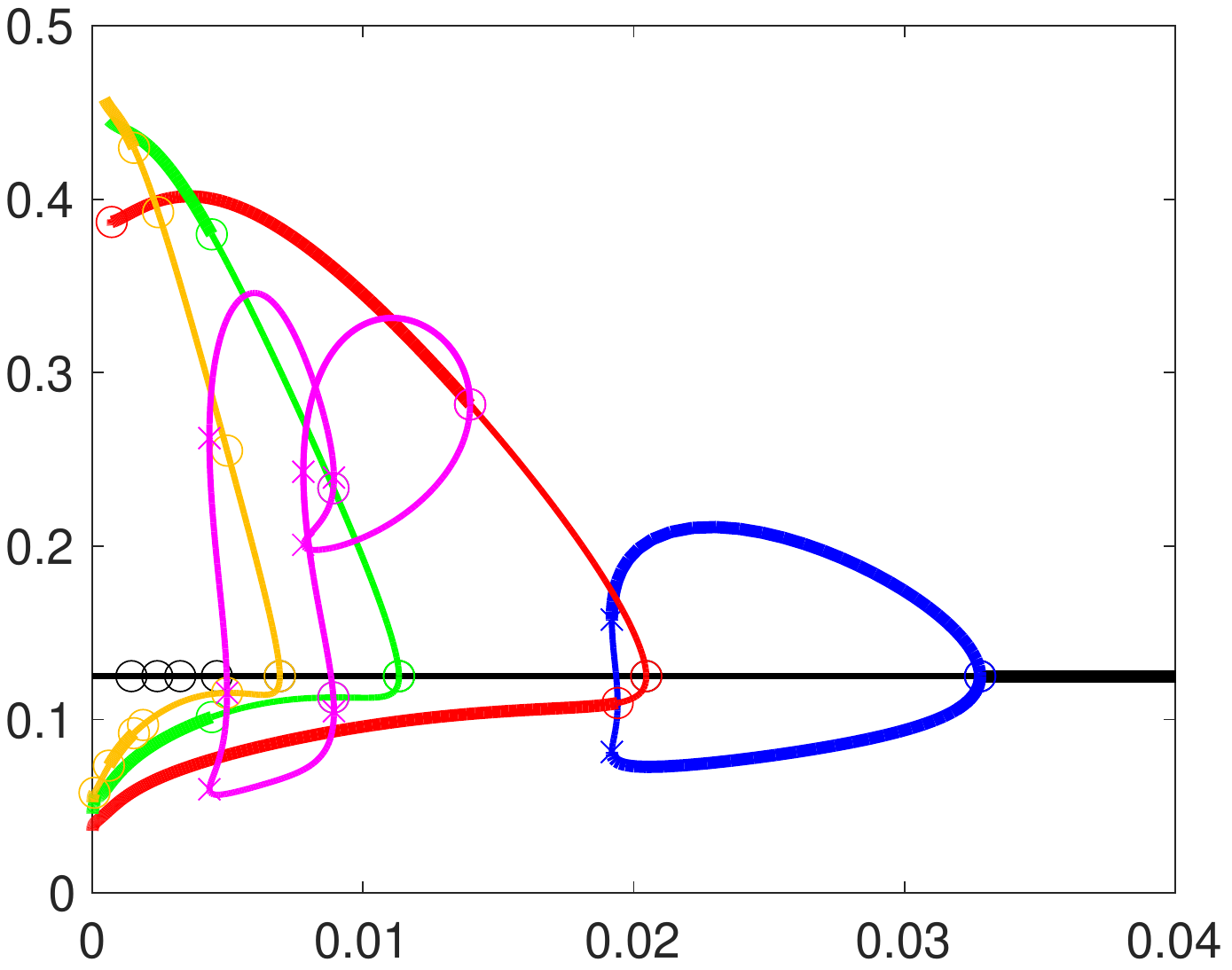}
\put(-10,68){$v(0)$}
\put(85,-2){$d$}
\end{overpic} 
}\hspace{0.1cm}
\subfloat[$d_{21}=0.01$\label{BD_d21_0p01}]{
\begin{overpic}[width=0.45\textwidth,tics=10]{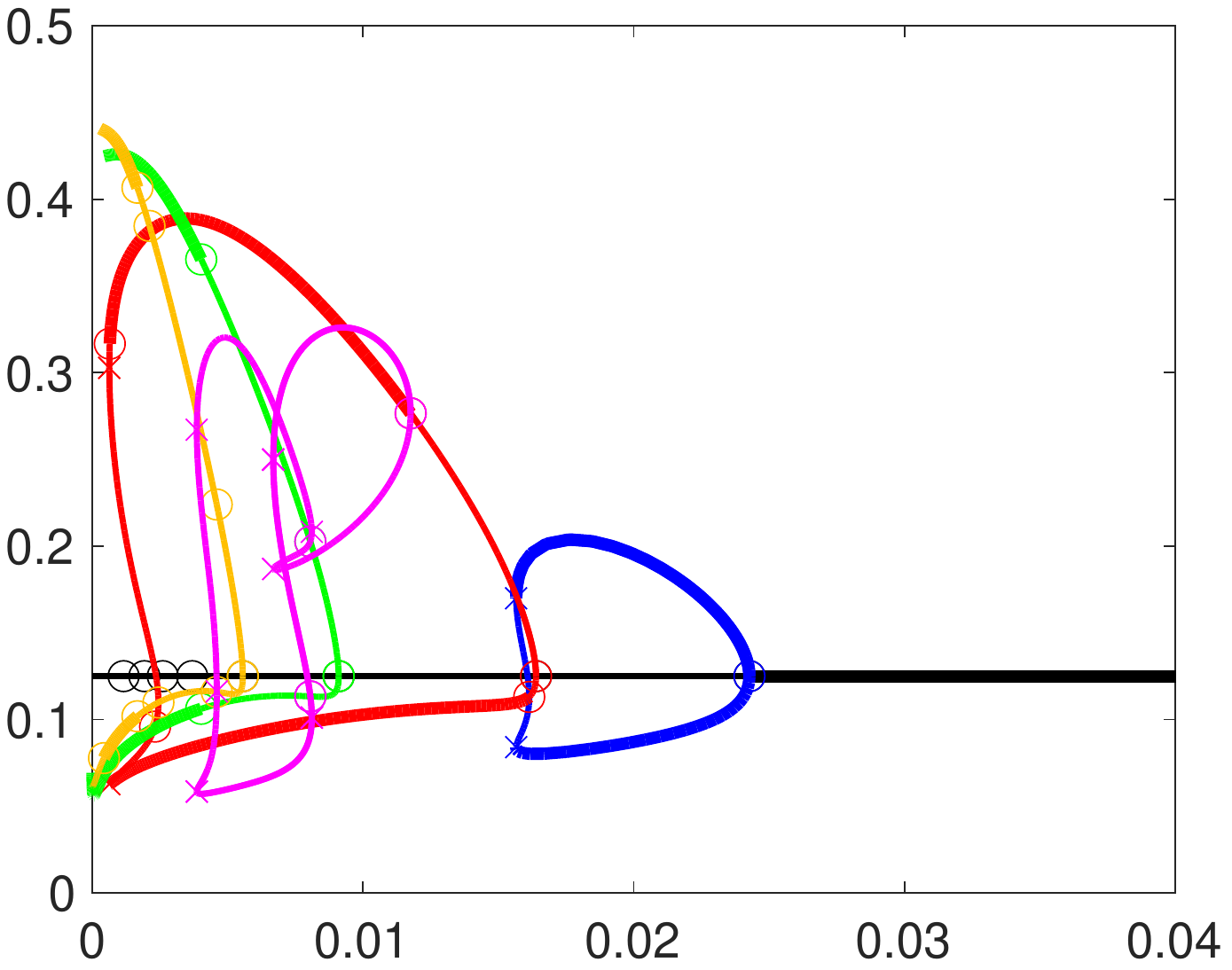}
\put(-10,68){$v(0)$}
\put(85,-2){$d$}
\end{overpic} 
}\\
\subfloat[$d_{21}=0.02$\label{BD_d21_0p02}]{
\begin{overpic}[width=0.45\textwidth,tics=10]{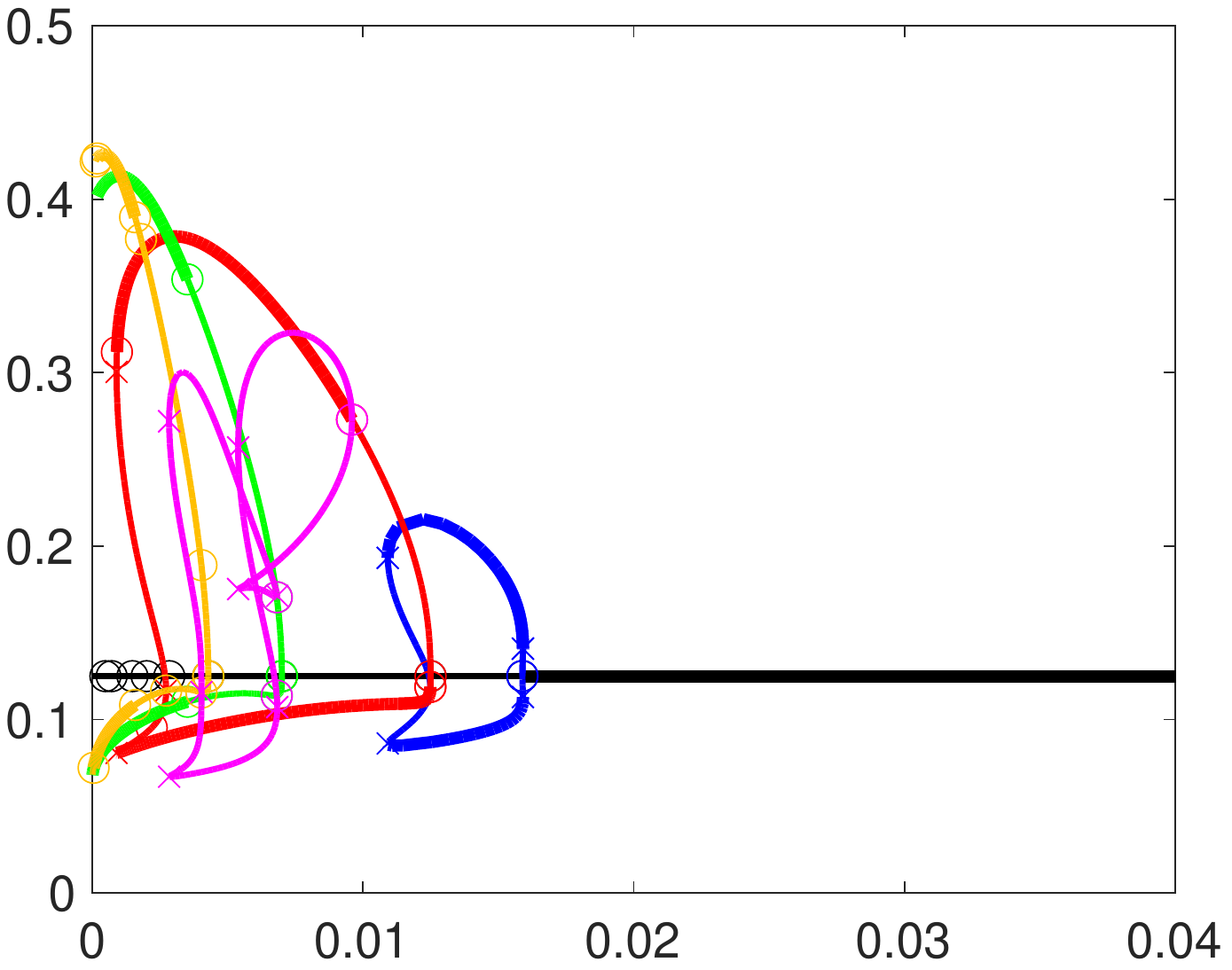}
\put(-10,68){$v(0)$}
\put(85,-2){$d$}
\end{overpic} 
}\hspace{0.1cm}
\subfloat[$d_{21}=0.025$\label{BD_d21_0p025}]{
\begin{overpic}[width=0.45\textwidth,tics=10]{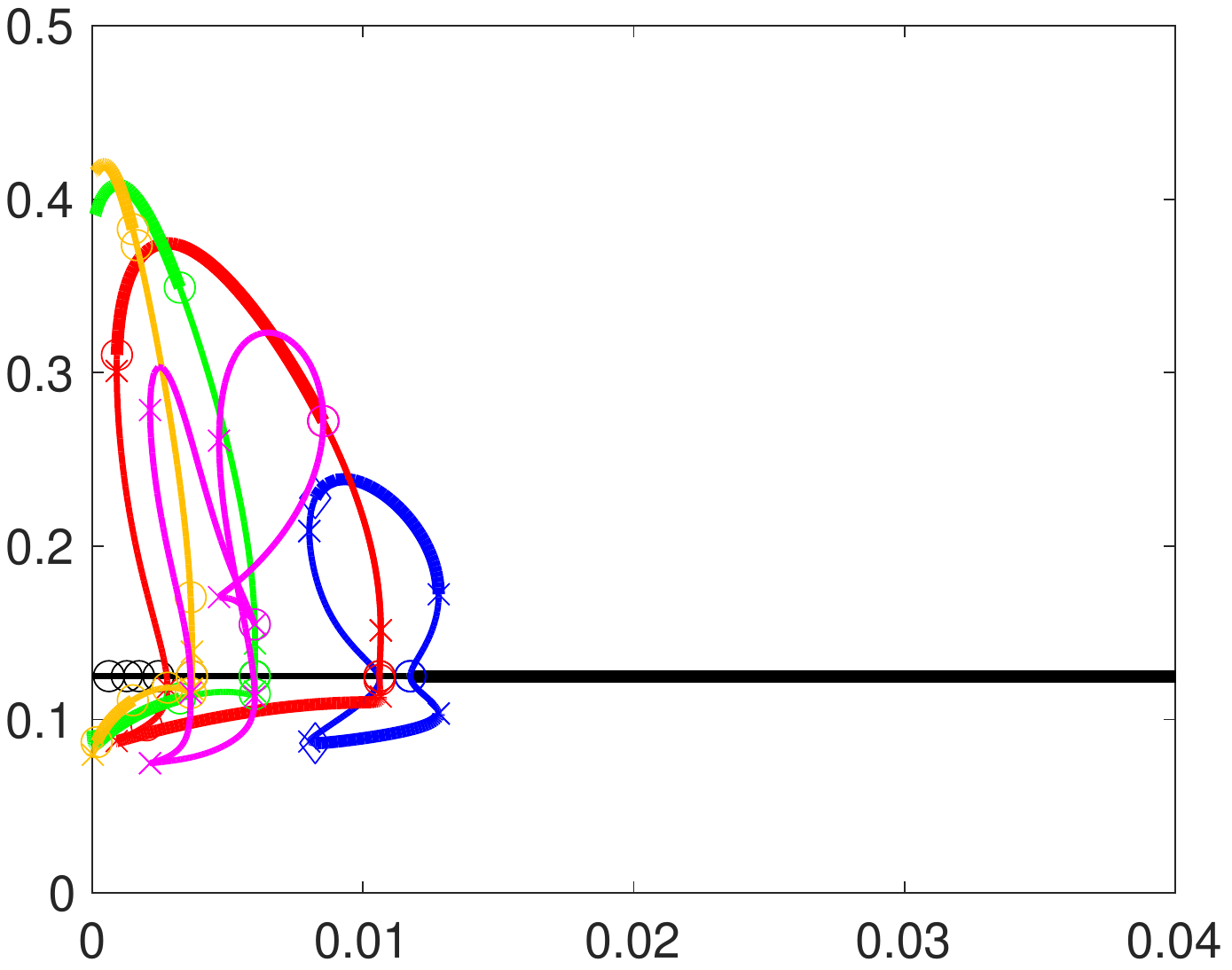}
\put(-10,68){$v(0)$}
\put(85,-2){$d$}
\end{overpic} 
}\\
\subfloat[$d_{21}=0.027$\label{BD_d21_0p027}]{
\begin{overpic}[width=0.45\textwidth,tics=10]{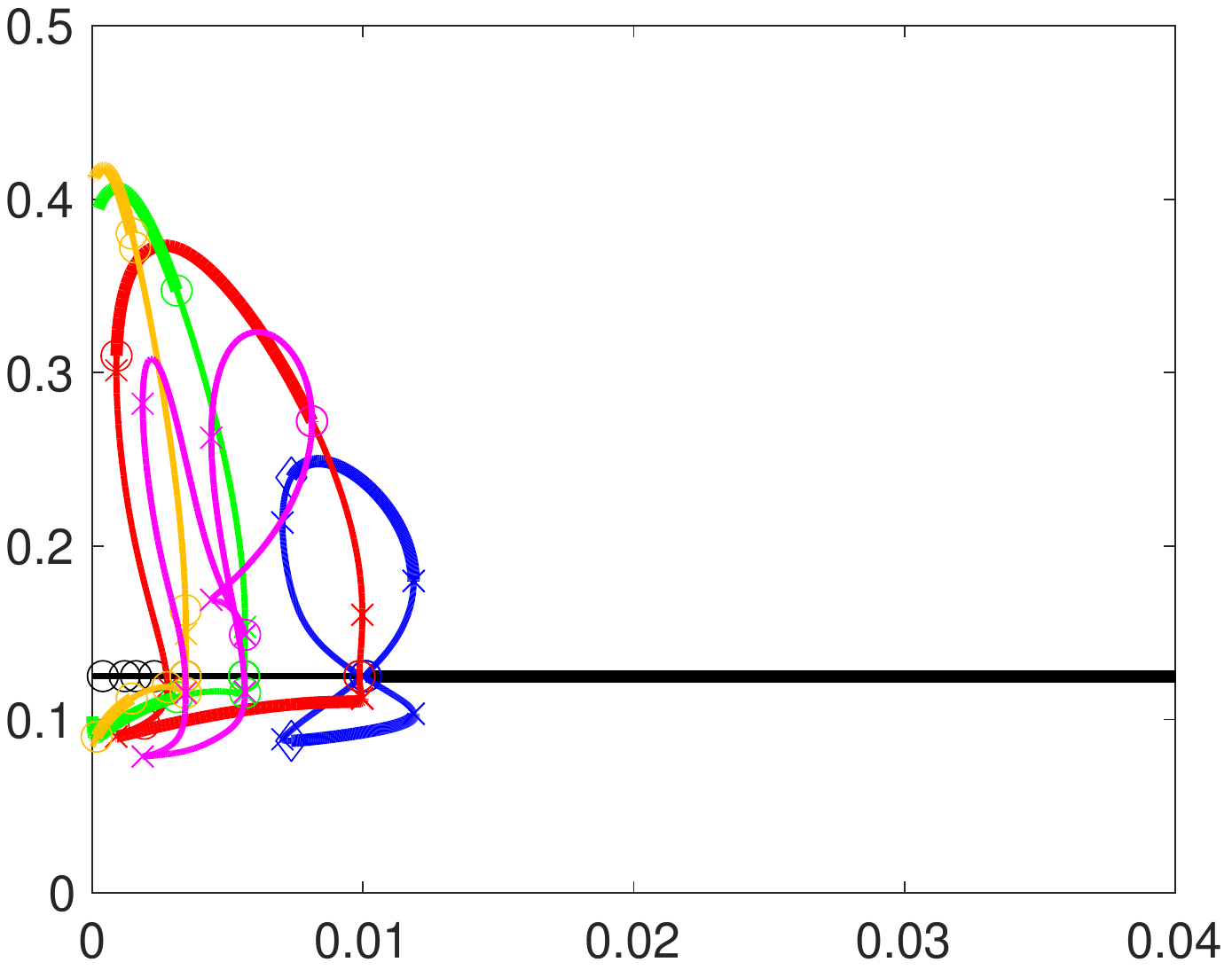}
\put(-10,68){$v(0)$}
\put(85,-2){$d$}
\end{overpic} 
}\hspace{0.1cm}
\subfloat[$d_{21}=0.03$\label{BD_d21_0p03}]{
\begin{overpic}[width=0.45\textwidth,tics=10]{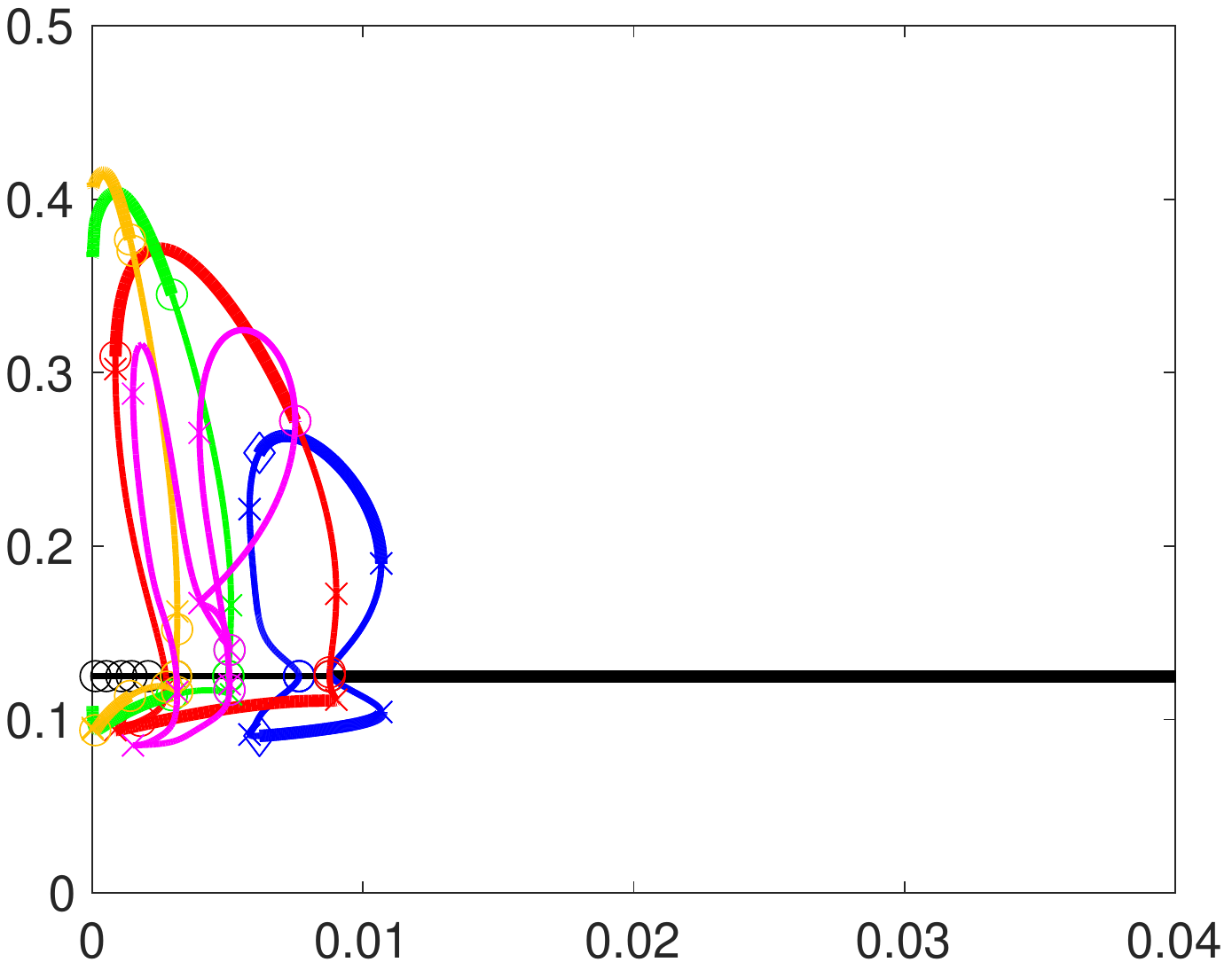}
\put(-10,68){$v(0)$}
\put(85,-2){$d$}
\end{overpic} 
}
\caption{Bifurcation diagrams for different values of the cross-diffusion coefficient $d_{21}$. The bifurcation parameter is the standard diffusion coefficient $d$, while on the $y$-axis we have $v(0)$. Thick/thin lines denotes stable/unstable stationary solutions. Circles/crosses/diamonds mark pitchfork/fold/Hopf bifurcations points. }\label{BD_d21}
\end{figure}
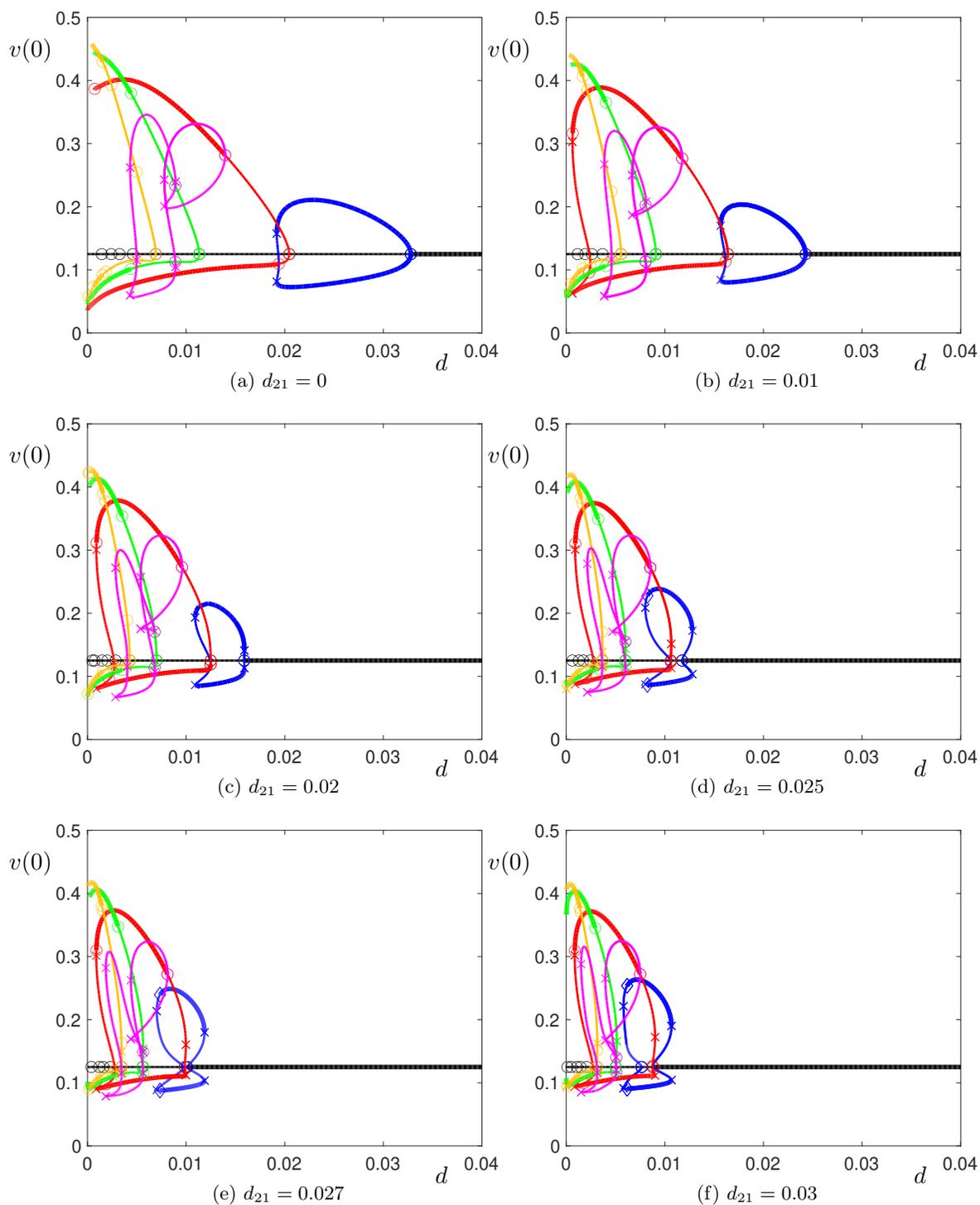

The interesting aspect is that the Hopf bifurcation point is located in a different position with respect to the doubly degenerate point than what observed in \cite{izuhara2018spatio}. Even though the increase of the cross-diffusion parameter $d_{21}$ shifts up the neutral stability curves, the appearance of the Hopf bifurcation point happens for $d_{12}>\hat{d}_{12}$ (``above'' the doubly degenerate point $(\hat{d},\hat{d}_{12})$, where the neutral stability curves associated to the 1- and 2-modes intersects), namely for parameter values for which the first two neutral stability curves have not yet switched their position (see the green dotted line in Figure \ref{nsc}).

\section{Weakly non-linear analysis}\label{sec_SL}

The goal of this section is to determine the type of pitchfork bifurcations on the homogeneous branch, depending on the parameter values. As observed, sub-critical pitchforks lead to a multi-stability scenario, where the system admits the stable homogeneous steady states together with one or more stable non-homogeneous ones. This can be done by deriving, through weakly nonlinear analysis, the Stuart--Landau equation for the amplitude of the patterns at the bifurcation point. We perform a weakly nonlinear analysis based on the method of multiple scales \cite{wollkind1994weakly}. The method is based on the fact that close to the bifurcation the amplitude of the pattern evolves on a slow temporal scale. Then new scaled coordinates are introduced and treated as separate variables in addition to the original variables \cite{wollkind1994weakly}. This has already been performed for system \eqref{cross} in \cite{gambino2012turing}, but for a different bifurcation parameter and with a different aim, namely to predict the amplitude of patterns. Here we consider the standard diffusion coefficient $d$ as bifurcation parameter, and the goal is to determine the type of pitchfork bifurcation and its dependence on the cross-diffusion coefficients.

The first important observation is that the linearised analysis gives information about the location of bifurcation points. To predict their type, non-linear terms must be included in the analysis.
The reaction cross-diffusion system \eqref{cross} (with $d_1=d_2=d$ and $d_{11}=d_{22}=0$) can be rewritten, separating linear and non-linear part, as
\begin{equation*}\label{cross-lin-nonlin}
\partial_t\bm{w}=\mathcal{L}^d\bm{w} +\dfrac{1}{2}Q_K(\bm{w},\bm{w})+\dfrac{1}{2}\Delta Q_D(\bm{w},\bm{w}),
\end{equation*}
where
\begin{equation*}
\bm{w}=\begin{pmatrix}u-u_*\\ v-v_*\end{pmatrix},
\end{equation*}
and the linear operator $\mathcal{L}^d$ is defined as
\begin{equation*}
\mathcal{L}^d(\bm{w})=K+D^d\Delta.
\end{equation*}
Note that the subscript indicates the dependence of the linear operator $\mathcal{L}$ and of the matrix $D$ on the bifurcation parameter $d$.
Being $\bm{x}=(x^u,x^v)$, $\bm{y}=(y^u,y^v),$ the bilinear operators that encodes the non-linear part are given by
\begin{equation*}
Q_K(\bm{x},\bm{y})=\begin{pmatrix}
-2a_1x^uy^u-b_1(x^uy^v+x^vy^u)\\
-2a_2x^vy^v-b_2(x^uy^v+x^vy^u)\end{pmatrix},\qquad
Q_D(\bm{x},\bm{y})=\begin{pmatrix}
d_{12}(x^uy^v+x^vy^u)\\
d_{21}(x^uy^v+x^vy^u)\end{pmatrix}.
\end{equation*}
Note that they are independent of the bifurcation parameter $d$.

We consider a bifurcation point at the value $d_c$ at which the mode relative to the eigenvalue $\lambda_k$ becomes unstable, which is generally indicated with $k_c$. Note that, once we focus on a particular eigenvalue $\lambda_k$, the value of $d_c$ is the positive root of \eqref{detM*d} given by \eqref{formula_d_c} (provided that condition \eqref{weak_comp} holds and $C_k$ in \eqref{ABC} is negative).

We now perform the weakly nonlinear expansion. Note that we follow \cite{gambino2012turing}, but we consider a different bifurcation parameter. Therefore the intermediate steps to derive the Stuart--Landau equation for the amplitude are the same (reported in the following for the reader's convenience), but the final coefficients will be different.

\begin{figure}
\centering
\subfloat[$L>0$\label{Lpos}]{
\begin{tikzpicture}
\draw (1.5,0) -- (3,0);
\draw[ultra thick] (3,0) -- (4.5,0);
\draw[ultra thick] (1.91,1.05) .. controls ++(-30:0.3) and ++(90:0.5) .. (3, 0)
.. controls ++(-90:0.5) and ++(30:0.2) .. (1.91,-1.05);
\draw[fill=white] (3,0) circle (0.9mm);
\draw (3.2,-0.4) node[] {$d_c$};
\draw (4.8,0) node[] {$d$};
\draw (1.2,0) node[] {\color{white}$d$};

\draw[gray!30!white] (3,0) circle (1.5cm);
\end{tikzpicture}
}
\hspace{2cm}
\subfloat[$L<0$\label{Lneg}]{
\begin{tikzpicture}
\draw (1.5,0) -- (3,0);
\draw[ultra thick] (3,0) -- (4.5,0);
\draw (4.09,1.05) .. controls ++(-150:0.3) and ++(90:0.5) .. (3, 0)
.. controls ++(-90:0.5) and ++(150:0.2) .. (4.09,-1.05);
\draw[fill=white] (3,0) circle (0.9mm);
\draw (2.8,-0.4) node[] {$d_c$};
\draw (4.8,0) node[] {$d$};
\draw (1.2,0) node[] {\color{white}$d$};
\draw[gray!30!white] (3,0) circle (1.5cm);
\end{tikzpicture}
}
\caption{Qualitative representation of the bifurcation structure close to the bifurcation point, predicted by the Stuart--Landau equation \eqref{SLeq}.}
\label{fig:locBD_L}
\end{figure}
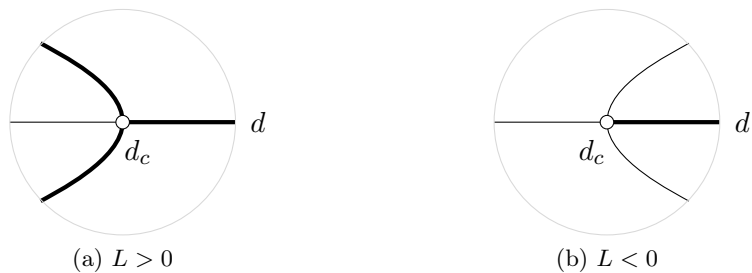

Close to the bifurcation, we can separate the fast time $t$ and slow time $T=\varepsilon^2 t$, where $\varepsilon$ is a small control parameter representing the dimensionless distance to the threshold. (see \cite{gambino2012turing} for more details). Following and therefore the time derivative decouples as $\partial t + \varepsilon^2 \partial T$.
Therefore, the time derivative decouples as $\partial t + \varepsilon^2 \partial T$, and we expand $d$ and $\bm{w}$ as
\begin{align*}
d&=d_c+\varepsilon^2d^{(2)}+\mathcal{O}(\varepsilon^4),\\
\bm{w}&=\varepsilon \bm{w} _1 + \varepsilon^2 \bm{w} _2 + \varepsilon^3 \bm{w}_3 +\mathcal{O}(\varepsilon^4),
\end{align*}
the linear operator
\begin{equation*}
\mathcal{L}^d=\mathcal{L}^{d_c}+\varepsilon^2d^{(2)} \begin{pmatrix} 1 & 0\\ 0& 1\end{pmatrix}\Delta +\mathcal{O}(\varepsilon^4),
\end{equation*}
and the bilinear operators (which are independent of $d$)
\begin{align*}
Q_K(\bm{w},\bm{w})&=\varepsilon^2 Q_K(\bm{w}_1,\bm{w}_1) +2 \varepsilon^3 Q_K(\bm{w}_1,\bm{w}_2) +\mathcal{O}(\varepsilon^4),\\
Q_D(\bm{w},\bm{w})&=\varepsilon^2 Q_D(\bm{w}_1,\bm{w}_1) +2 \varepsilon^3 Q_D(\bm{w}_1,\bm{w}_2) +\mathcal{O}(\varepsilon^4).
\end{align*}
We substitute the expansions in \eqref{cross-lin-nonlin}, obtaining the following equation
\begin{equation*}
\begin{split}
\varepsilon^3 \partial_T \bm{w}_1=
\varepsilon &\mathcal{L}^{d_c}\bm{w}_1+
\varepsilon^2 \mathcal{L}^{d_c}\bm{w}_2+ \varepsilon^3 \mathcal{L}^{d_c}\bm{w}_3
+\varepsilon^3 d^{(2)} \begin{pmatrix} 1 & 0\\ 0& 1\end{pmatrix}\Delta \bm{w}_1\\
&+\dfrac 1 2 \varepsilon^2 Q_K(\bm{w}_1,\bm{w}_1) +\varepsilon^3 Q_K(\bm{w}_1,\bm{w}_2)\\
&+\dfrac 1 2 \varepsilon^2 \Delta Q_D(\bm{w}_1,\bm{w}_1) +\varepsilon^3 \Delta Q_K(\bm{w}_1 ,\bm{w}_2)
+\mathcal{O}(\varepsilon^4),
\end{split}
\end{equation*}
from which we get the following equations for $\bm{w}_1,\,\bm{w}_2$ and $\bm{w}_3$, collecting all the terms corresponding to the same order in $\varepsilon$:
\begin{align}
\boxed{\varepsilon^1} \qquad & \mathcal{L}^{d_c}\bm{w}_1=0,\label{order1}\\
\boxed{\varepsilon^2} \qquad & \mathcal{L}^{d_c}\bm{w}_2=F,\label{order2}\\
\boxed{\varepsilon^3} \qquad & \mathcal{L}^{d_c}\bm{w}_3=G,\label{order3}
\end{align}
where
\begin{align}
\label{eq:F}
F&=-\dfrac 1 2 \varepsilon^2 Q_K(\bm{w}_1,\bm{w}_1)-\dfrac 1 2 \varepsilon^2 \Delta Q_D(\bm{w}_1,\bm{w}_1),\\
\label{eq:G}
G&=\partial_T \bm{w}_1-\varepsilon^3 d^{(2)} \begin{pmatrix} 1 & 0\\ 0& 1\end{pmatrix}\Delta \bm{w}_1-Q_K(\bm{w}_1,\bm{w}_2)-\Delta Q_K(\bm{w}_1,\bm{w}_2).
\end{align}
A solution to the first (linear) equation \eqref{order1} satysfying homogeneous Neumann boundary condition is
\begin{equation}\label{w_1}
\bm{w}_1=A(T)\bm{\rho} \cos(k_cx)
\end{equation}
with $\bm{\rho} \in \Ker(K-k_c^2D^{d_c})$. Then (using the fact that $\det(K-k_c^2D^{d_c})=0$) we find
\begin{equation*}
\rho=\begin{pmatrix} 1\\ M\end{pmatrix}, \qquad M=\dfrac{K_{21}-k_c^2D^{d_c}_{21}}{D^{d_c}_{22}k_c^2-K_{22}}.
\end{equation*}
Note that, from \eqref{w_1}, we have that the leading term $\bm{w}_1$ is the product of a slowly varying amplitude and the basic pattern.
We solve now the second equation \eqref{order2}, where $F$ given in \eqref{eq:F} can be rewritten as
\begin{equation*}
F=-\dfrac 1 4 A^2 \sum_{i=1,2}\mathcal{M}_i(\bm{\rho},\bm{\rho})\cos(ik_cx), \qquad
\mathcal{M}_i(\rho,\rho):= Q_K(\bm{\rho},\bm{\rho})-i^2k_c^2Q_D(\bm{\rho},\bm{\rho}).
\end{equation*}
From the Fredholm alternative, equation \eqref{order2} has a solution if and only if
\begin{equation*}
\forall \bm{\psi}\in \Ker((K-k_cD^{d_c})^T) \quad \langle F,\bm{\psi} \rangle=0.
\end{equation*}
We have that
\begin{equation*}
\bm{\psi}=\begin{pmatrix}1\\M_*\end{pmatrix}\cos(k_cx), \qquad M_*=\dfrac{K_{12}-k_c^2D^{d_c}_{12}}{D^{d_c}_{22}k_c^2-K_{22}},
\end{equation*}
again thanks to the fact that $\det(K-k_c^2D^{d_c})=0$, and then it is easy to check that the solvability condition is verified (in $L_2(0,1/k_c)$) without additional conditions.
Then we have that a solution is
\begin{equation*}
\bm{w}_2=A^2\sum_{i=0,2}\bm{w}_{2i}\cos(ik_cx),
\end{equation*}
where $\bm{w}_{2i}, \, i=0,2$ solve
\begin{equation*}
\left(K-i^2k_c^2 D^{d_c} \right)\bm{w}_{2i}=-\dfrac 1 4 \mathcal{M}(\bm{\rho},\bm{\rho}), \quad i=0,2.
\end{equation*}
Regarding the last equation \eqref{order3}, we rewrite $G$ given in \eqref{eq:G} as
\begin{equation*}
G=\left( \dfrac{dA}{dT}\bm{\rho}+AG_1^{(1)}+A^3G_1^{(3)}\right)\cos(k_c x)
+A^3G_3 cos(3k_c x),
\end{equation*}
with
\begin{align*}
G_1^{(1)}&=d^{(2)}k_c^2\bm{\rho},\\
G_1^{(3)}&=-\mathcal{M}_1(\bm{\rho},\bm{w}_{20})-\dfrac{1}{2}\mathcal{M}_1(\bm{\rho},\bm{w}_{22})\\
G_3&=-\dfrac{1}{2}\mathcal{M}_3(\bm{\rho},\bm{w}_{22}).
\end{align*}
The solvability condition $\langle G,\bm{\psi}\rangle=0$ translates into the Stuart--Landau equation for the amplitude $A(T)$
\begin{equation}\label{SLeq}
\dfrac{dA}{dT}=\sigma A- L A^3,
\end{equation}
where the coefficients $\sigma$ and $L$ are given by
\begin{equation*}
\sigma:=-\dfrac{\langle G_1^{(1)},\begin{pmatrix} 1 \\ M_*\end{pmatrix}\rangle}{\langle \bm{\rho},\begin{pmatrix} 1 \\ M_*\end{pmatrix}\rangle}=-d^{(2)}k_c^2, \qquad L:=\dfrac{\langle G_1^{(3)},\begin{pmatrix} 1 \\ M_*\end{pmatrix}\rangle}{\langle\bm{\rho},\begin{pmatrix} 1 \\ M_*\end{pmatrix}\rangle},
\end{equation*}
Note that, when $L<0$, this expansion is not able to predict the amplitude of the patterns and higher-order terms must be included in the weakly nonlinear analysis. However, we are only here interested in the type of bifurcation, and equation \eqref{SLeq} it is enough to this task. In detail, when $L>0$, at $d_c$ the transition occurs via super-critical bifurcation, while when $L<0$, the bifurcation is sub-critical. The type of bifurcation (super- or sub-critical) predicted by the Stuart--Landau equation \eqref{SLeq} depending on the sign of the parameter $L$ is reported in Figure \ref{fig:locBD_L}.

The expression for $\sigma$ is very clear and simple, on the contrary, the expression for $L$ is intricate and it is not possible to analytically obtain information about its dependence on the cross-diffusion coefficients. We have now to evaluate it numerically. The Matlab scripts generating the following Figures are available at \cite{GitHubSL}.

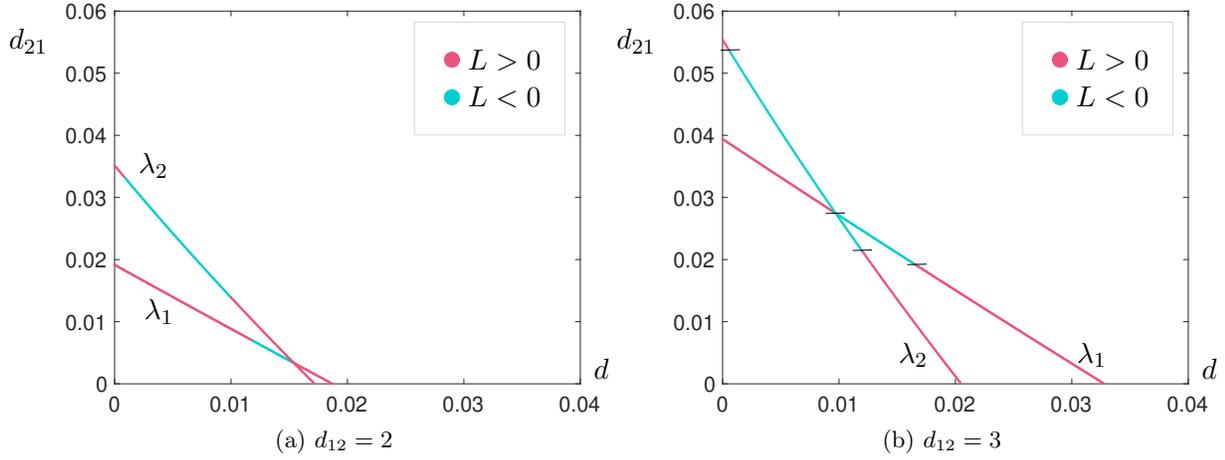
\begin{figure}
\centering
\subfloat[$d_{12}=2$\label{sgnL_d_d21_with_d12_2}]{
\begin{overpic}[width=0.45\textwidth,tics=10]{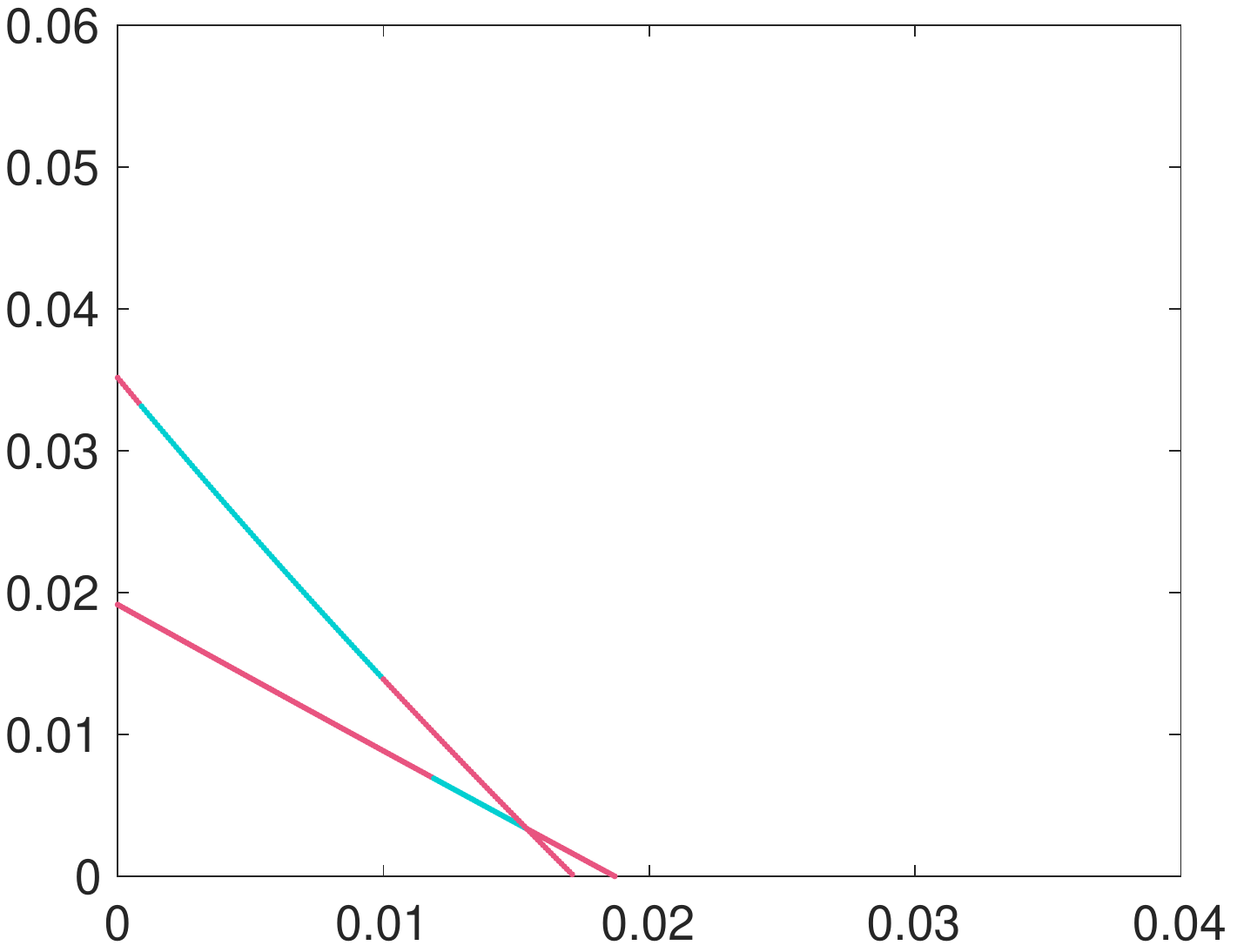}
\put(-10,68){$d_{21}$}
\put(98,7){$d$}
\put(15,18){$\lambda_1$}
\put(14,45){$\lambda_2$}
\put(63,52){
\begin{tikzpicture}
\draw[Lpos,fill=Lpos] (0,0) circle (0.1cm);
\draw (0.7,0) node {$L>0$};
\draw[Lneg,fill=Lneg] (0,-0.5) circle (0.1cm);
\draw (0.7,-0.5) node {$L<0$};
\draw[gray!30!white] (-0.5,0.5) rectangle (1.5,-1);
\end{tikzpicture}
}
\end{overpic}
}\hspace{0.5cm}
\subfloat[$d_{12}=3$\label{sgnL_d_d21_with_d12_3}]{
\begin{overpic}[width=0.45\textwidth,tics=10]{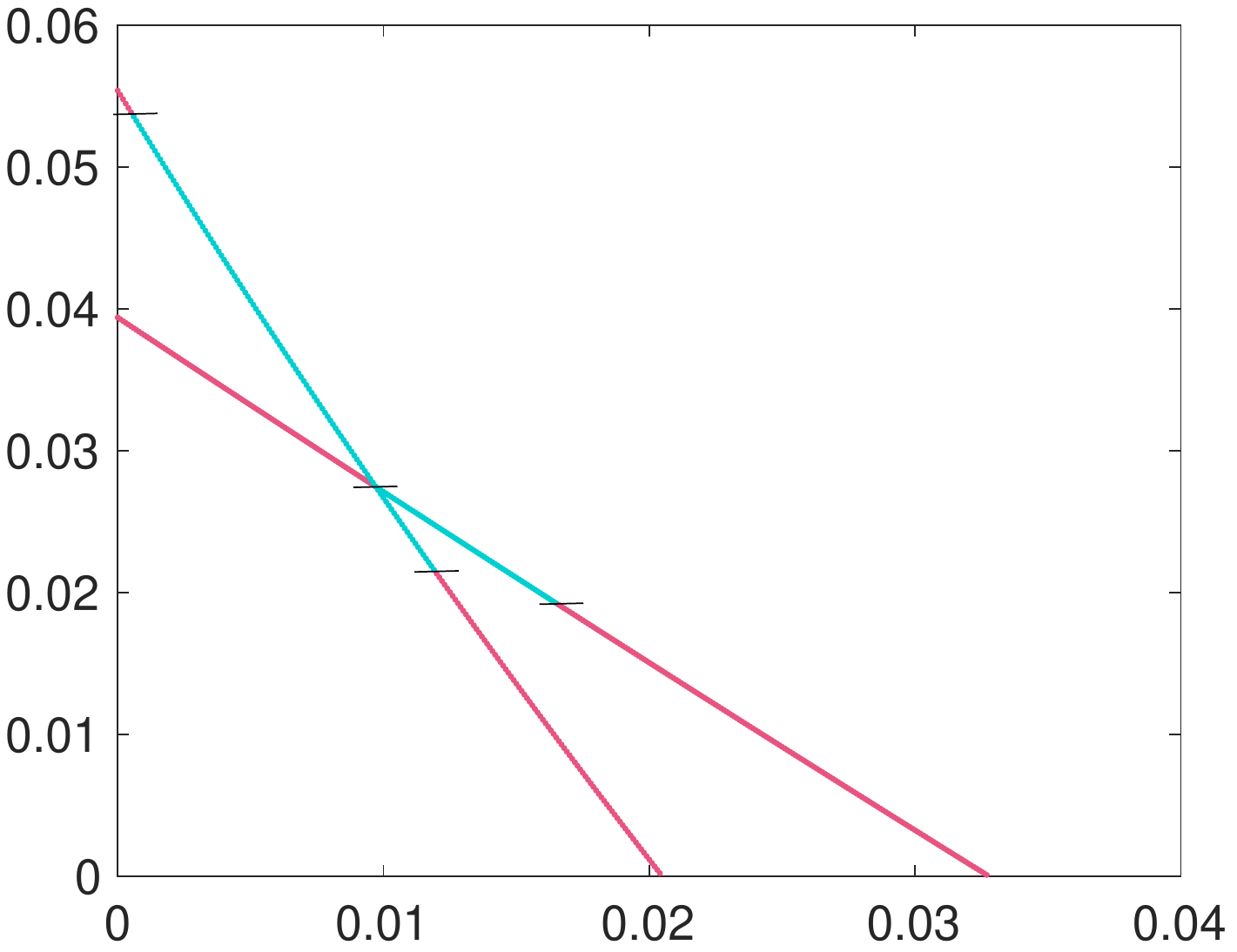}
\put(-10,68){$d_{21}$}
\put(98,7){$d$}
\put(75,10){$\lambda_1$}
\put(42,10){$\lambda_2$}
\put(63,52){
\begin{tikzpicture}
\draw[Lpos,fill=Lpos] (0,0) circle (0.1cm);
\draw (0.7,0) node {$L>0$};
\draw[Lneg,fill=Lneg] (0,-0.5) circle (0.1cm);
\draw (0.7,-0.5) node {$L<0$};
\draw[gray!30!white] (-0.5,0.5) rectangle (1.5,-1);
\end{tikzpicture}
}
\end{overpic}
}
\caption{Sign of $L$ along the neutral stability curves of $\lambda_1$ and $\lambda_2$ in the $(d,d_{21})$-plane, with $d_{12}=2$ (left) and $d_{12}=3$ (right). The remain parameter values are listed in Table \ref{tab:param}.}\label{sgnL_d_d21_lambda12}
\end{figure}

The results relevant to the first two eigenvalues of the Laplacian along the neutral stability curves in the $(d,d_{21})$-plane for two values of $d_{12}$ (to better explain the bifurcation diagrams in Figure \ref{BD_d21}) are shown in Figure \ref{sgnL_d_d21_lambda12}. On the neutral stability curves $d_{21}=d_{21}(d_c,d_{12},\lambda_k)$, the sign of $L$ is indicated with different colours, while small black horizontal lines mark the changes of sign along the curves. With $d_{12}=2$, we observe a change of sign of $L(\lambda_1)$ at the doubly degenerate point, while grater values of $d_{21}$ triggers the appearance of multi-stable solutions (homogenous and 2-mode non-homogenous). With $d_{12}=3$, as observed in the bifurcation diagrams in Section \ref{sec_bif}, when $d_{21}$ increases, the branch corresponding to $\lambda_1$ from super-critical becomes sub-critical, and we have a multi-stability region where the stable homogeneous steady state coexists with stable non-homogeneous (1-mode) steady states. By increasing further $d_{21}$, the same happen to the branch corresponding to $\lambda_2$. Finally, when the neutral stability curves cross at $(\hat{d},\hat{d}_{21})$, the first bifurcation point is related to the 2-mode, and it turns out to be a sub-critical pitchfork, while the sign of $L$ of the bifurcation point corresponding to $\lambda_1$ becomes positive again. This corresponds to a multi-stability region where the stable homogeneous steady state coexists with stable non-homogeneous (1-mode and 2-mode) steady states. These different cases can be observed also in Figure \ref{BD_d21}.

We also note that for $d_{21}$ the appearance of the Hopf bifurcation point on the $\lambda_1$ branch happens exactly when the coefficient $L(\lambda_2)$ changes sign. However, this happens only for the particular value $d_{12}=3$, meaning that the two events are not related.

\begin{figure}
\centering
\subfloat[$\lambda_1$\label{sgnL_d12_d21_lambda1}]{
\begin{overpic}[width=0.45\textwidth,tics=10]{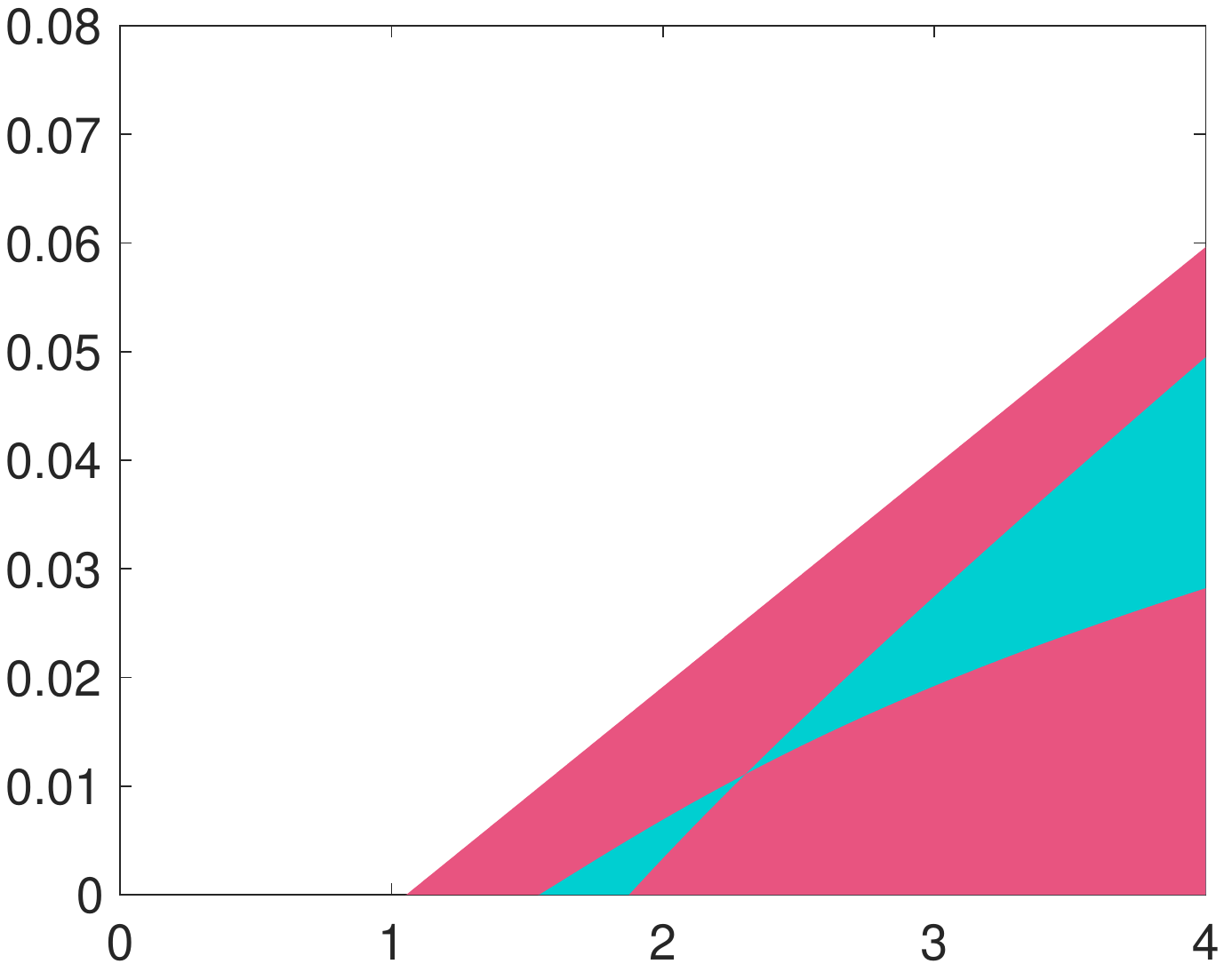}
\put(-10,68){$d_{21}$}
\put(85,-2){$d_{12}$}
\put(11,53){
\begin{tikzpicture}
\draw[Lpos,fill=Lpos] (0,0) circle (0.1cm);
\draw (0.7,0) node {$L>0$};
\draw[Lneg,fill=Lneg] (0,-0.5) circle (0.1cm);
\draw (0.7,-0.5) node {$L<0$};
\draw[gray!30!white] (-0.5,0.5) rectangle (1.5,-1);
\end{tikzpicture}
}
\end{overpic}
}\hspace{0.5cm}
\subfloat[$\lambda_2$\label{sgnL_d12_d21_lambda2}]{
\begin{overpic}[width=0.45\textwidth,tics=10]{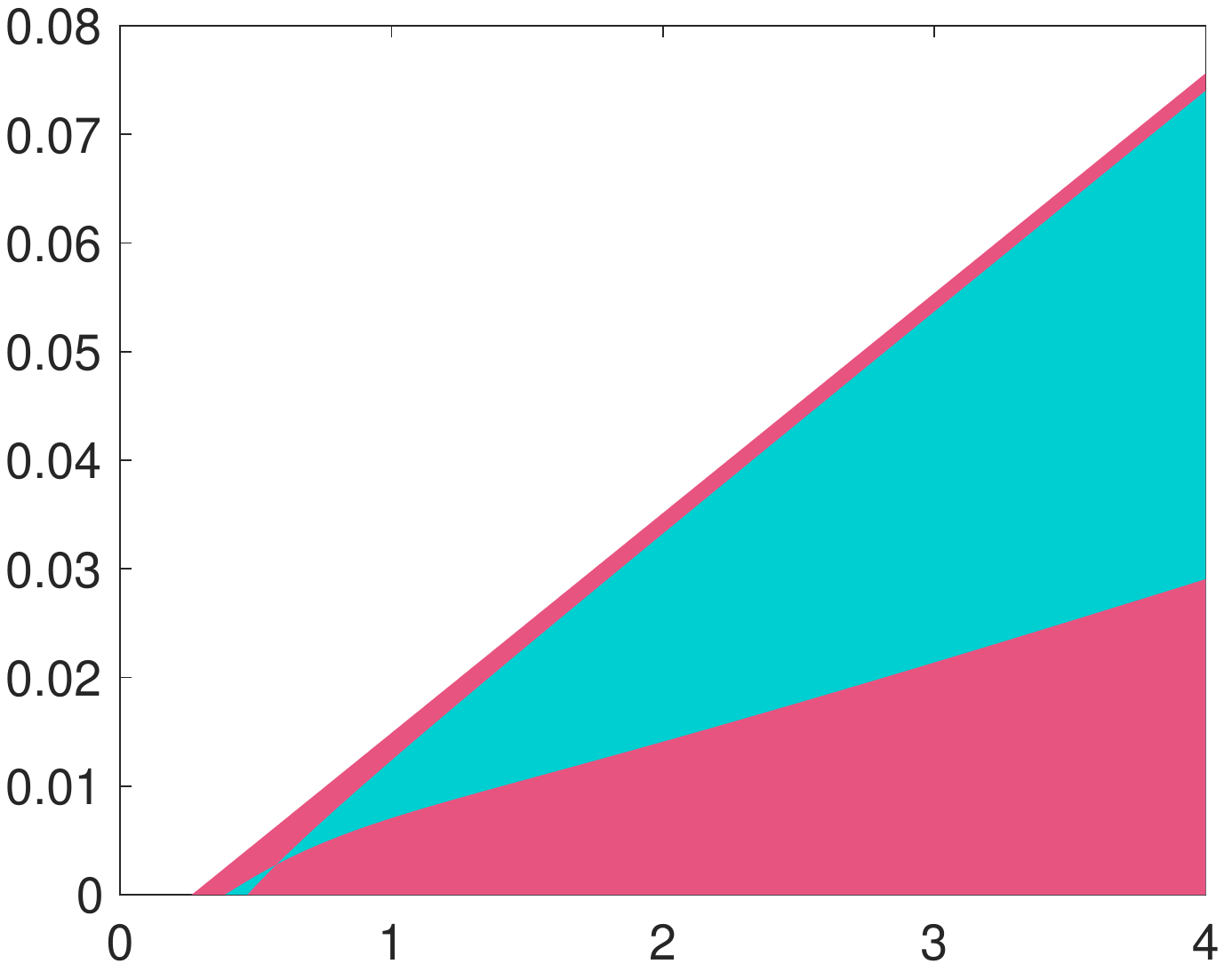}
\put(-10,68){$d_{21}$}
\put(85,-2){$d_{12}$}
\put(11,53){
\begin{tikzpicture}
\draw[Lpos,fill=Lpos] (0,0) circle (0.1cm);
\draw (0.7,0) node {$L>0$};
\draw[Lneg,fill=Lneg] (0,-0.5) circle (0.1cm);
\draw (0.7,-0.5) node {$L<0$};
\draw[gray!30!white] (-0.5,0.5) rectangle (1.5,-1);
\end{tikzpicture}
}
\end{overpic}
}
\caption{Sign of $L$ for $\lambda_1$ (left) and $\lambda_2$ (right) in the $(d_{12},d_{21})$-plane.}\label{sgnL_d12_d21}
\end{figure}
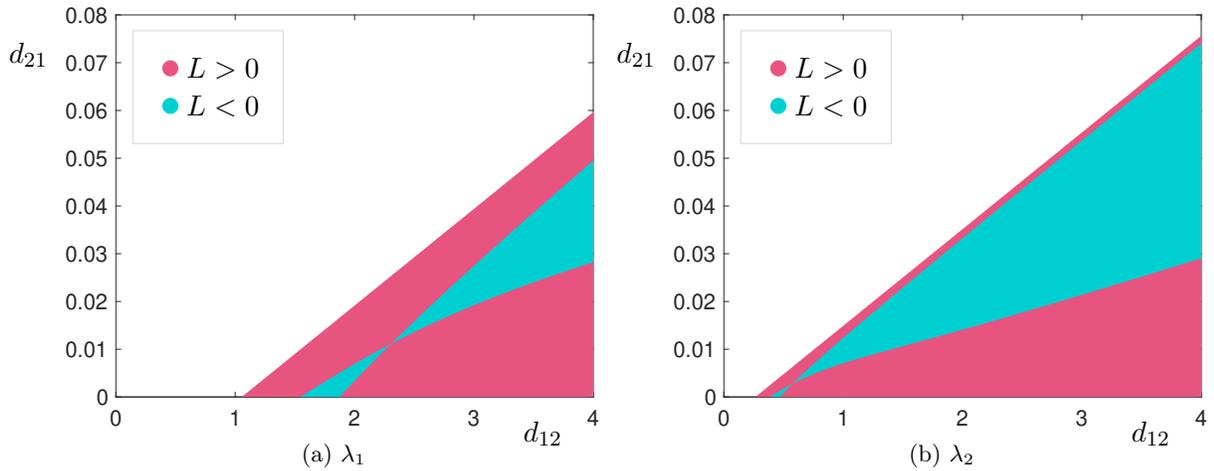

The bifurcation structure is not peculiar of the particular parameter set considered, and the sign of $L$ at the bifurcation points can be investigated for different values of the cross-diffusion coefficients. This allows to detect the regions where multi-stability is present. In Figure \ref{sgnL_d12_d21}, the sign of $L$ in the $(d_{12},d_{21})$-plane is shown. Note that $d_{12}$ determine a threshold value for $d_{21}$ to the existence of the bifurcation point related to an eigenvalue $\lambda_k$ \cite{breden2021influence}.

\section{Hopf bifurcation points}\label{sec_Hopf}
In this section, we investigate the presence of Hopf bifurcation points in the bifurcation structure of system~\eqref{cross} (remember that we neglect self-diffusion, namely~$d_{11}=d_{22}=0$). In~\cite{izuhara2018spatio}, the existence of spatially non-constant time-periodic solutions has been proven in a rigorous way for the triangular case, appearing from the doubly degenerate point $(\hat{d},\hat{d}_{12})$ (intersection of two neutral stability curves). Note that, when the additional cross-diffusion is turned on, the coordinates of the doubly degenerate point changes for increasing values of~$d_{21}$.

Also in this context, we mainly focus on the influence of the cross-diffusion coefficient $d_{21}$ and on how the bifurcation structure modifies when it increases. In Section~\ref{sec_bif}, we have already observed that, for the fixed value~$d_{12}=3$, the increasing of the other cross-diffusion coefficient leads to the appearance of Hopf bifurcation points, that are not present in the triangular case. Looking at the neutral stability curves in Figure~\ref{neutral_stab_curves}, we see that we are in the vicinity of the doubly degenerate point, but with~$d_{12}>\hat{d}_{12}$. This suggests that another effect may have occurred. In particular, one of the differences between this case and the one reported in \cite{izuhara2018spatio} (namely the~$d_{12}=1.7$ and~$d_{21}=0$), consists in the mutual position and type of pitchfork bifurcations.

For this reason, we investigated in more detail the sign of L for the 1- and 2-mode along the neutral stability curves. The results in the~$(d,d_{12})$-plane for different values of $d_{21}$ are shown in Figure~\ref{nsc_signL_d_d12_d21}. In the triangular case (Figure~\ref{nsc_signL_d_d12_d21_0}), the neutral stability curves intersect at~$(\hat{d},\hat{d}_{12})$: along the $\lambda_1$-curve, the sign of $L(\lambda_1)$ changes at the intersection. In detail, $L(\lambda_1)$ is negative for $d_{12}<\hat{d}_{12}$ and $d<\hat{d}$. From \cite{izuhara2018spatio}, we know that for~$d_{12}<\hat{d}_{12}$ (close to the doubly degenerate point) the bifurcation diagram presents a Hopf bifurcation point at $d_H$ and stable time-periodic spatial patterns arises for $d<d_H$. Increasing the value of $d_{21}$, the region with $L(\lambda_1)<0$ reduces (Figure~\ref{nsc_signL_d_d12_d21_0p01}) and disappears (Figure~\ref{nsc_signL_d_d12_d21_0p01105}). In particular, $L(\lambda_1)$ is negative for $d_{12}>\hat{d}_{12}$ and $d>\hat{d}$ (Figure~\ref{nsc_signL_d_d12_d21_0p012}). Numerical investigations shows that in this case (close to the doubly degenerate point) the bifurcation diagram presents a Hopf bifurcation point at $d_H$ with both~$d_{12}<\hat{d}_{12}$ and~$d_{12}<\hat{d}_{12}$ and stable time-periodic spatial patterns arises for $d>d_H$. Finally, for grater values of $d_{21}$, also $L(\lambda_2)$ changes sign (Figures~\ref{nsc_signL_d_d12_d21_0p018},~\ref{nsc_signL_d_d12_d21_0p025}), without qualitative modifications in the bifurcation structure.

\begin{figure}
\centering
\subfloat[\label{beforeHH}]{
\begin{tikzpicture}
\draw (1.5,0) -- (4.5,0);
\draw[ultra thick,Lpos] (4.05,1.05) .. controls ++(60:-0.3) and ++(30:0.5) .. (3, 0);
\draw[ultra thick,Lneg] (3, 0) .. controls ++(30:-0.5) and ++(15:0.5) .. (1.65,-0.65);
\draw[ultra thick,gray] (3.45,1.44) .. controls ++(80:-0.6) and ++(70:0.5) .. (3, 0)
.. controls ++(70:-0.5) and ++(30:0.5) .. (2.1,-1.2);
\draw[fill=white] (3,0) circle (0.9mm);
\draw[gray!30!white] (3,0) circle (1.5cm);
\draw[thick,yellow] (4.525,0) arc (0:-180:1.525);
\end{tikzpicture}
}
\hspace{1cm}
\subfloat[\label{atHH}]{
\begin{tikzpicture}
\draw (1.5,0) -- (4.5,0);
\draw[ultra thick,Lpos] (4.05,1.05) .. controls ++(60:-0.3) and ++(30:0.5) .. (3, 0);
\draw[ultra thick,Lpos] (3, 0) .. controls ++(30:-0.5) and ++(15:0.5) .. (1.65,-0.65);
\draw[ultra thick,,gray] (3.45,1.44) .. controls ++(80:-0.6) and ++(70:0.5) .. (3, 0)
.. controls ++(70:-0.5) and ++(30:0.5) .. (2.1,-1.2);
\draw[fill=white] (3,0) circle (0.9mm);
\draw[gray!30!white] (3,0) circle (1.5cm);
\draw[ultra thick, dotted, white] (3,0) circle (1.54cm);
\end{tikzpicture}
}
\hspace{1cm}
\subfloat[\label{afterHH}]{
\begin{tikzpicture}
\draw (1.5,0) -- (4.5,0);
\draw[ultra thick,Lneg] (4.05,1.05) .. controls ++(60:-0.3) and ++(30:0.5) .. (3, 0);
\draw[ultra thick,Lpos] (3, 0) .. controls ++(30:-0.5) and ++(15:0.5) .. (1.65,-0.65);
\draw[ultra thick,gray] (3.45,1.44) .. controls ++(80:-0.6) and ++(70:0.5) .. (3, 0)
.. controls ++(70:-0.5) and ++(30:0.5) .. (2.1,-1.2);
\draw[fill=white] (3,0) circle (0.9mm);
\draw[gray!30!white] (3,0) circle (1.5cm);
\draw[ultra thick, dotted, yellow] (3,0) circle (1.525cm);
\end{tikzpicture}
}
\caption{Qualitative representation of sign of $L$ ({\color{Lpos} $\bullet$} $L>0$, {\color{Lneg} $\bullet$} $L<0$), predicted by the Stuart--Landau equation \eqref{SLeq}, along the neutral stability curves close to the doubly degenerate point. The region in which time-periodic spatial patter may appear is marked in yellow (solid line denotes stable solutions, dotted line unstable ones).}
\label{fig:loc_L_H}
\end{figure}
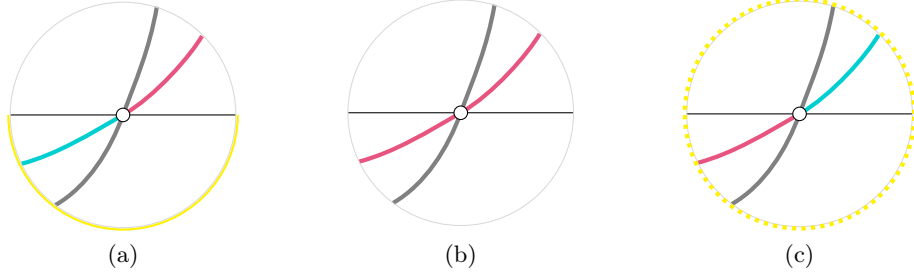

A qualitative representation of the changes happening along the neutral stability curves for increasing values of $d_{21}$ are sketched in Figure~\ref{fig:loc_L_H}.

\begin{figure}
\centering
\subfloat[$d_{21}=0$ \label{nsc_signL_d_d12_d21_0}]{
\begin{overpic}[width=0.45\textwidth,tics=10]{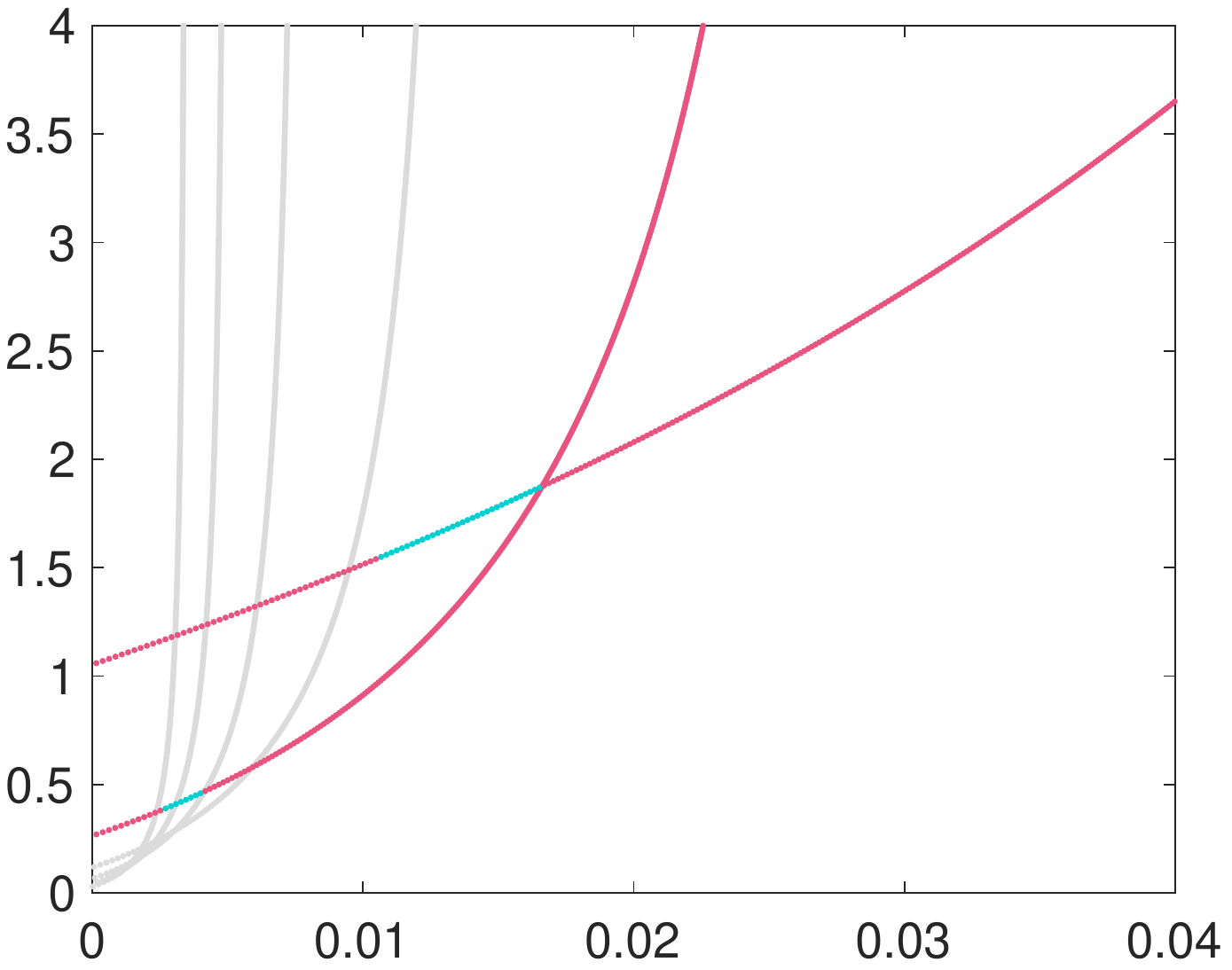}
\put(-8,60){$d_{12}$}
\put(97,7){$d$}
\put(85,70){$\lambda_1$}
\put(47,70){$\lambda_2$}
\end{overpic}
}\hspace{0.1cm}
\subfloat[$d_{21}=0.01$\label{nsc_signL_d_d12_d21_0p01}]{
\begin{overpic}[width=0.45\textwidth,tics=10]{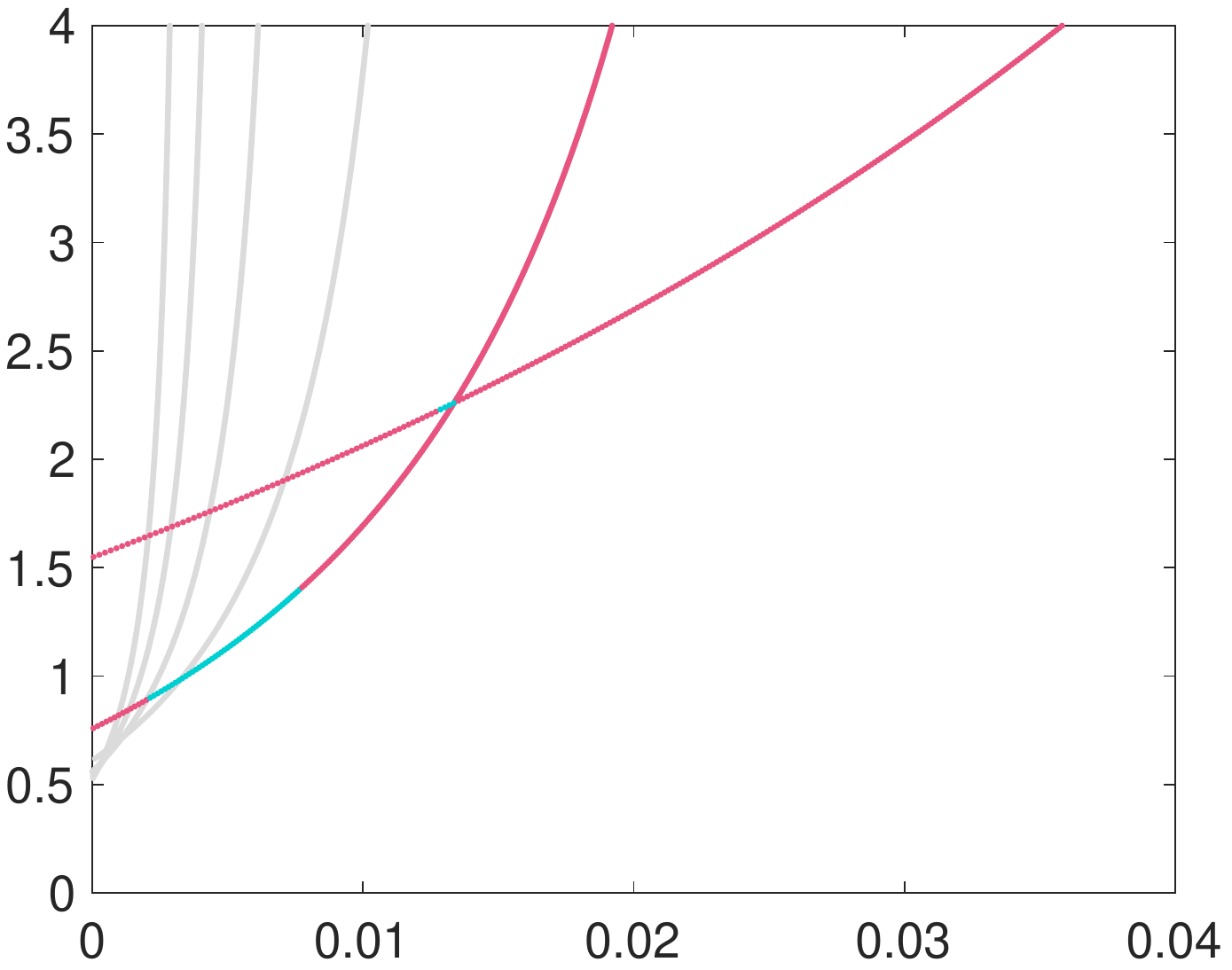}
\put(-8,60){$d_{12}$}
\put(97,7){$d$}
\put(85,70){$\lambda_1$}
\put(50,70){$\lambda_2$}
\end{overpic}
}\\
\subfloat[$d_{21}=0.01105$ \label{nsc_signL_d_d12_d21_0p01105}]{
\begin{overpic}[width=0.45\textwidth,tics=10]{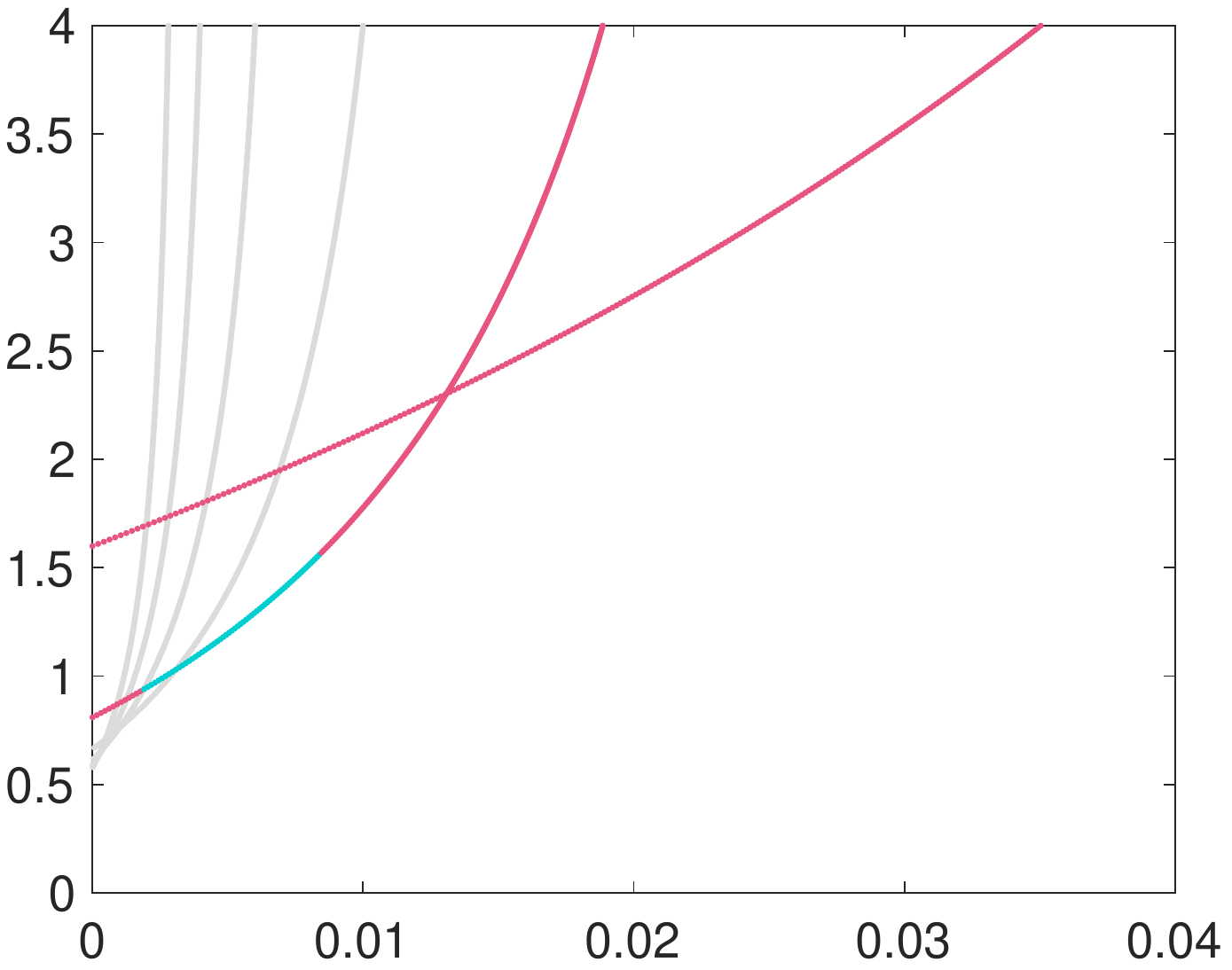}
\put(-8,60){$d_{12}$}
\put(97,7){$d$}
\put(85,70){$\lambda_1$}
\put(50,70){$\lambda_2$}
\end{overpic}
}\hspace{0.1cm}
\subfloat[$d_{21}=0.012$\label{nsc_signL_d_d12_d21_0p012}]{
\begin{overpic}[width=0.45\textwidth,tics=10]{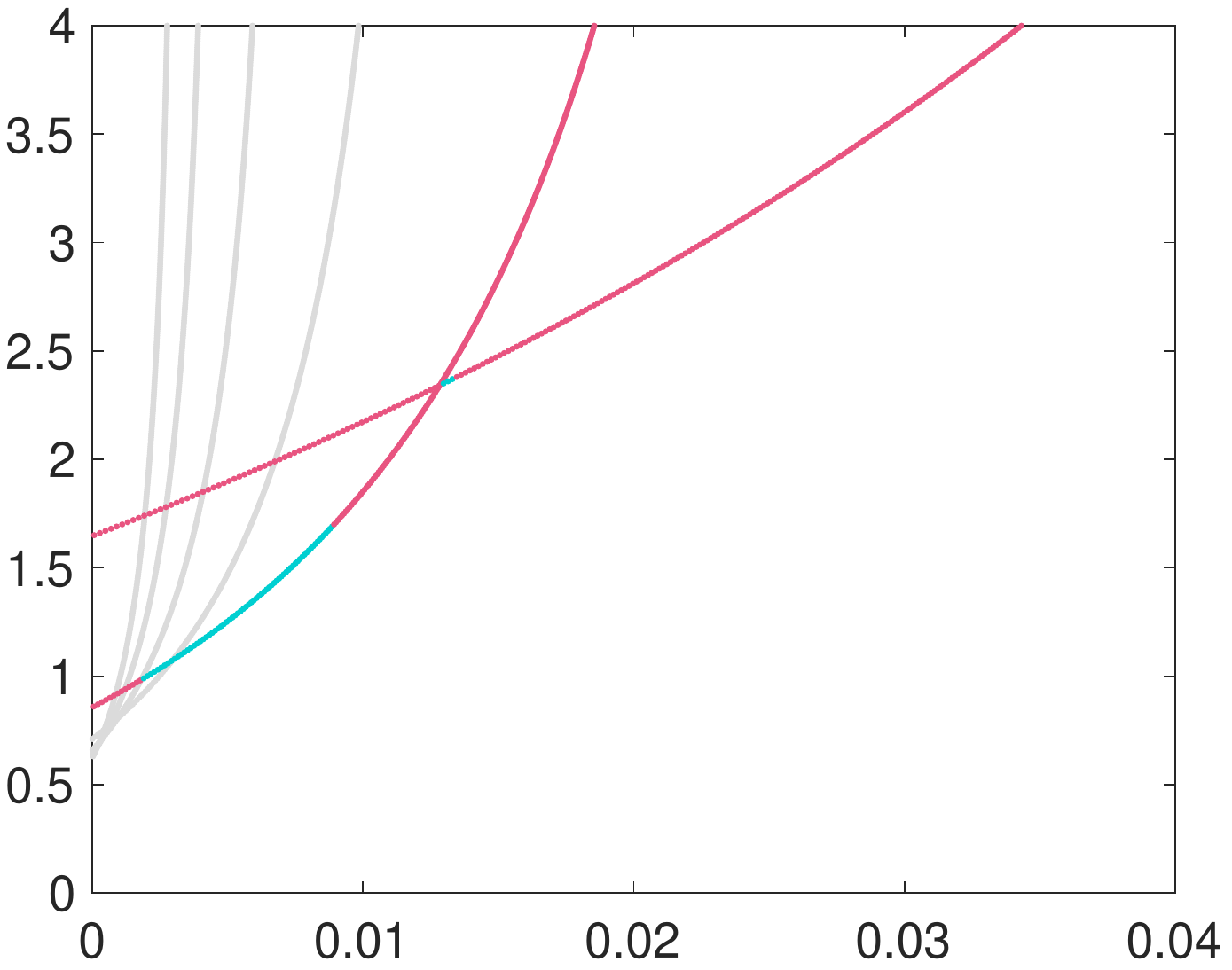}
\put(-8,60){$d_{12}$}
\put(97,7){$d$}
\put(85,70){$\lambda_1$}
\put(50,70){$\lambda_2$}
\end{overpic}
}\\
\subfloat[$d_{21}=0.018$ \label{nsc_signL_d_d12_d21_0p018}]{
\begin{overpic}[width=0.45\textwidth,tics=10]{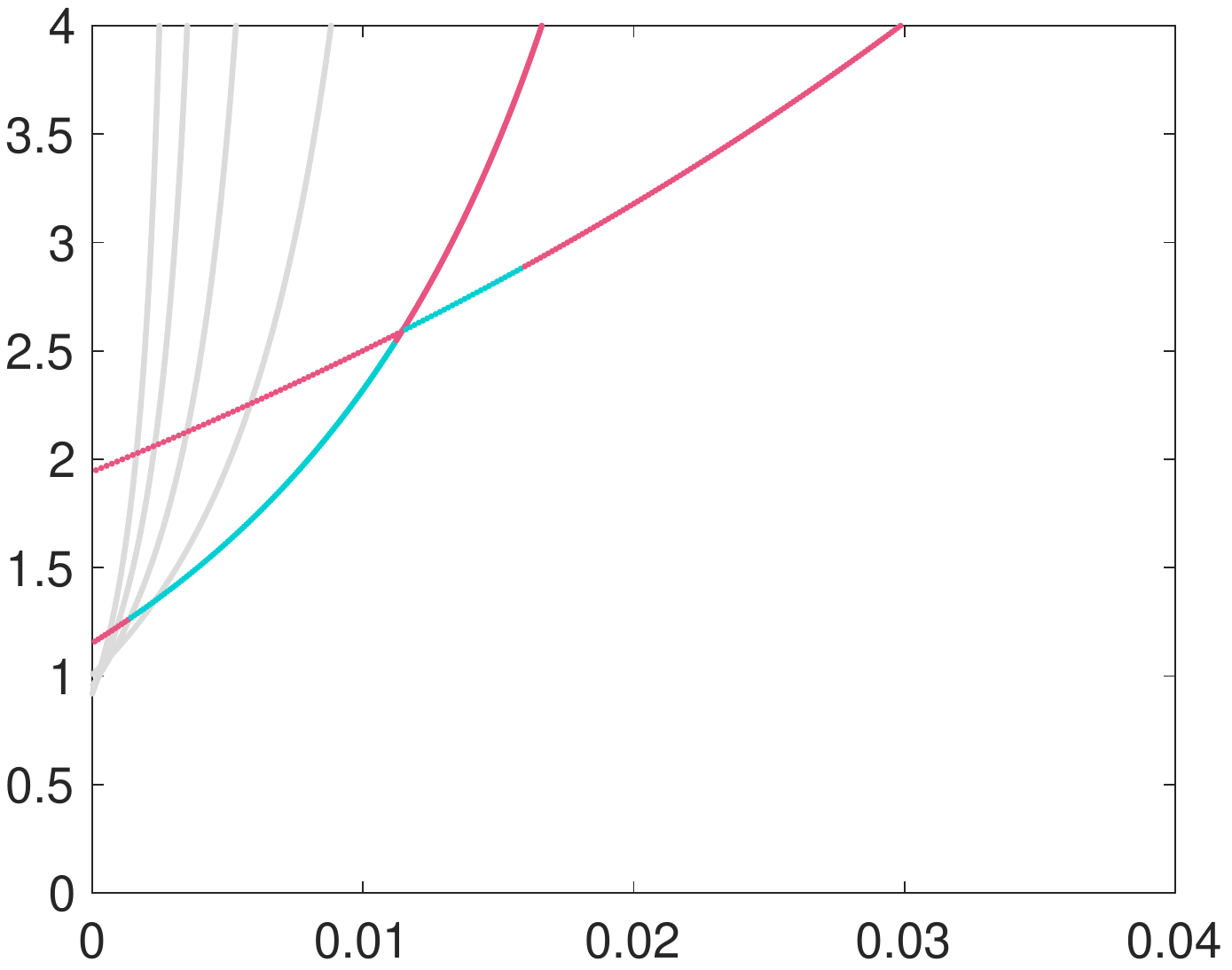}
\put(-8,60){$d_{12}$}
\put(97,7){$d$}
\put(70,70){$\lambda_1$}
\put(45,70){$\lambda_2$}
\end{overpic}
}\hspace{0.1cm}
\subfloat[$d_{21}=0.025$\label{nsc_signL_d_d12_d21_0p025}]{
\begin{overpic}[width=0.45\textwidth,tics=10]{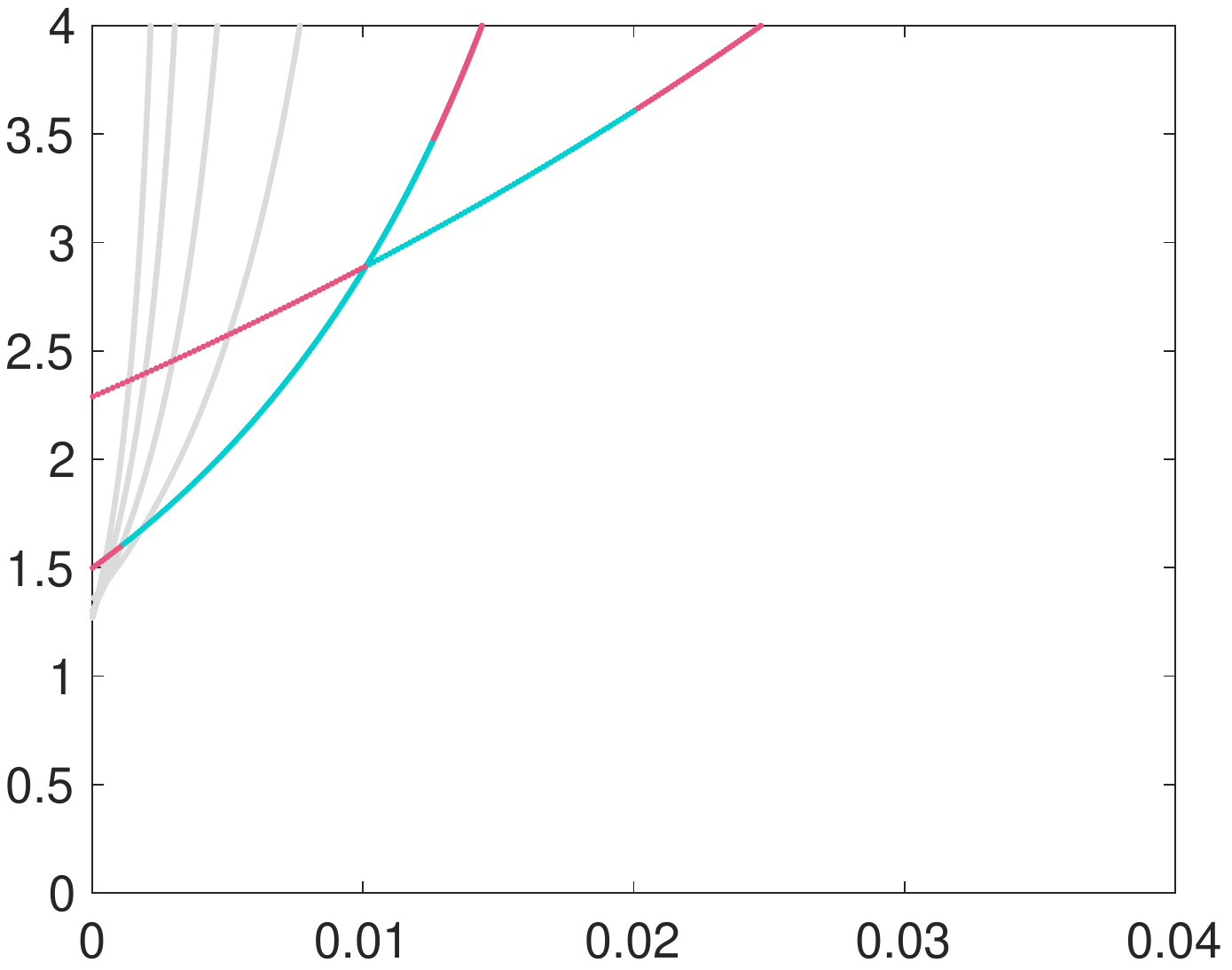}
\put(-8,60){$d_{12}$}
\put(97,7){$d$}
\put(65,70){$\lambda_1$}
\put(40,70){$\lambda_2$}
\end{overpic}
}
\caption{Sign of the coefficient $L$ along the neutral stability curves for $\lambda_k,\,k=1,\dots,6$ for different values of the cross-diffusion coefficient $d_{21}$. Colours appears on the curves related to the first two modes, ({\color{Lpos} $\bullet$} $L>0$, {\color{Lneg} $\bullet$} $L<0$), while the other modes are marked in gray.
}\label{nsc_signL_d_d12_d21}
\end{figure}
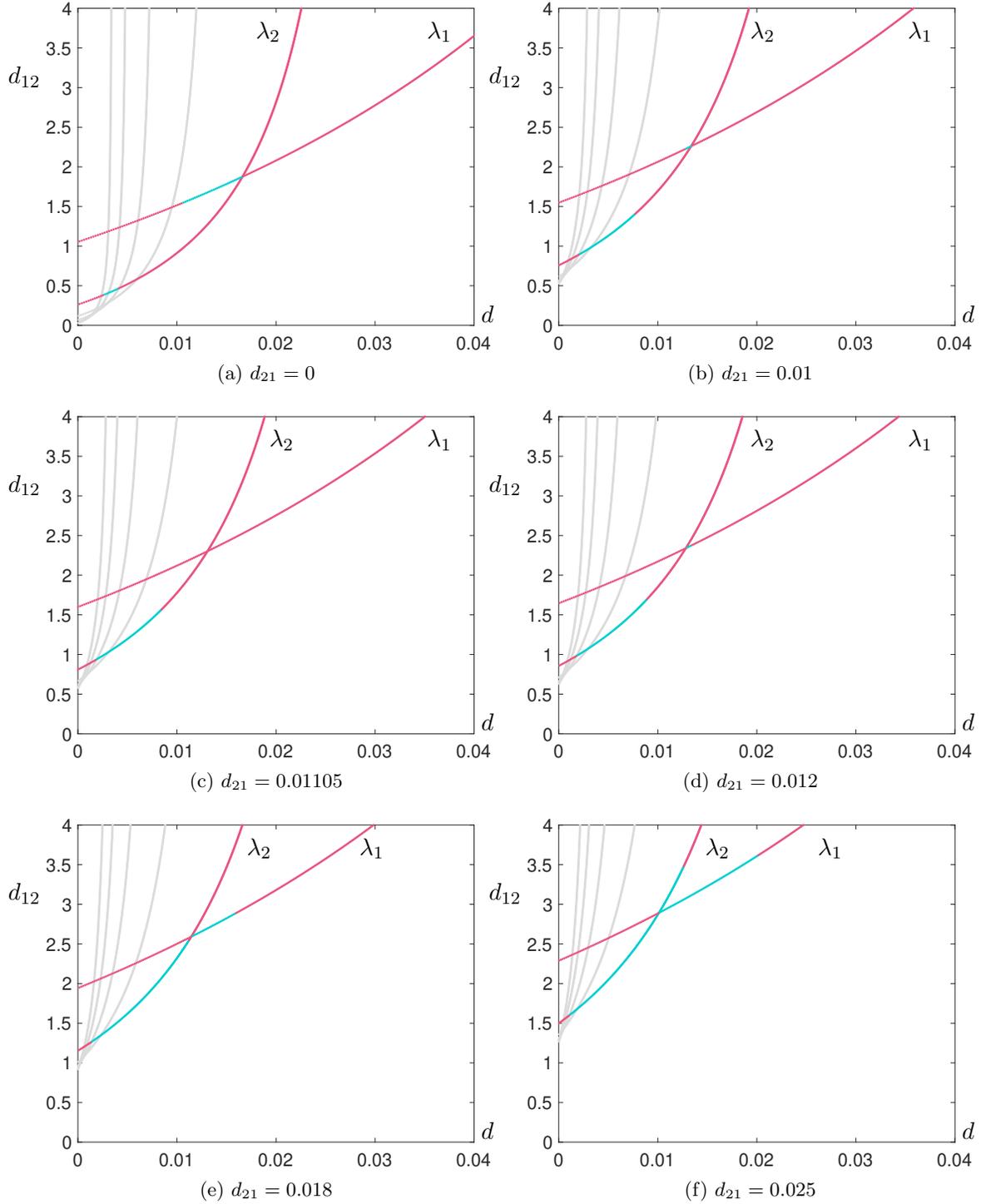

Finally, we want to generalise the method proposed in \cite{izuhara2018spatio}, where the center manifold reduction has been applied to study the dynamics around the doubly degenerate point in the triangular case.

Following the technique presented for the triangular case, the cross-diffusion system is transformed into a infinite dimensional dynamical system in the Fourier modes
\begin{equation*}
\begin{pmatrix}\dot{x}_k\\ \dot{y}_k \end{pmatrix} =M_k \begin{pmatrix}u_k\\v_k\end{pmatrix} + \begin{pmatrix}F_k\\ G_k\end{pmatrix},\quad k\in \mathbb{N}_0,
\end{equation*}
with
\begin{align*}
F_k&=-a_1\sum_{k_1+k_2=k}u_{k_1}u_{k_2}-\sum_{k_1+k_2=k}(d_{12}\lambda_k+b_1)u_{k_1}v_{k_2},\\
G_k&=-a_2\sum_{k_1+k_2=k}v_{k_1}v_{k_2}-\sum_{k_1+k_2=k}(d_{21}\lambda_k+b_2)v_{k_1}v_{k_2}.
\end{align*}
In the full cross-diffusion system the matrices $T_k$ that diagonalise $M_k,\, k=1,2$ given in \ref{matrixM_k} depends on $d_{21}$ and are given by
\begin{equation*}
T_k:=\begin{pmatrix}T^k_{11}&T^k_{12}\\T^k_{21} & T^k_{22}\end{pmatrix}
=\begin{pmatrix}M^k_{12}&M^k_{11}\\-M^k_{11} & M^k_{21}\end{pmatrix}, \quad k=1,2.
\end{equation*}
This is the only (and natural) modification needed to the method proposed for the triangular case. Setting
\begin{equation*}
\begin{pmatrix}x_k\\y_k\end{pmatrix}
=T^{-1}_k\begin{pmatrix}u_k\\v_k\end{pmatrix}, \quad k=1,2,
\end{equation*}
the dynamics around the doubly degenerate point is described by the following infinite dimensional dynamical system
\begin{align}
&\dot{x}_k=\dfrac{1}{\det T_k}\left( T^k_{22}\tilde{F}_k-T^k_{12}\tilde{G}_k\right), &&k=1,2,\nonumber\\
&\dot{y}_k=(\tr M_k)y_k+\dfrac{1}{\det T_k}\left(-T^k_{21}\tilde{F}_k+T^k_{11}\tilde{G}_k\right), && k=1,2,\label{cross_ddp}\\
&\begin{pmatrix}\dot{x}_k\\ \dot{y}_k \end{pmatrix} =M_k \begin{pmatrix}u_k\\v_k\end{pmatrix} + \begin{pmatrix}F_k\\ G_k\end{pmatrix},&& k\in \mathbb{N}_0\setminus\{1,2\},\nonumber
\end{align}
where $\tilde{F}_k$ and $\tilde{G}_k$ denote the nonlinear terms $F_k$ and $G_k$ depending on $(x_k,y_k)^T$.
Sufficiently close to the doubly degenerate point $(\hat{d},\hat{d}_{12})$ the dynamics of \eqref{cross_ddp} on the center manifold is topologically equivalent to
\begin{align*}
\dot{x}_1&=\mu_1x_1+A_1x_1x_2+(A_2x_1^2+A_3x_2^2)x_1+\mathcal{O}(|(x_1,x_2)|^4),\\
\dot{x}_2&=\mu_2x_1+B_1x_1^2+(B_2x_1^2+B_3x_2^2)x_2+\mathcal{O}(|(x_1,x_2)|^4),
\end{align*}
where the coefficients $A_k,\,B_k \in \mathbb{R}$ are explicitly determined (see \cite{izuhara2018spatio} for more details) and they depend on the parameters of both the reaction and the diffusion part (in particular they depend on $d_{21}$). In particular, we have
\begin{equation*}
A_1=\dfrac{1}{\det T_1}\left(T^1_{22}f^1_1-T^1_{12}g^1_1\right), \qquad B_1=\dfrac{1}{\det T_2}\left(T^2_{22}f^1_2-T^1_{12}g^1_2\right),
\end{equation*}
and
\begin{align*}
f^1_1&=-2a_1T^1_{11}T^2_{11}-(\hat{d}_{12}\lambda_1+b_1)(T^1_{11}T^2_{21}+T^2_{11}T^1{21}),& f^1_2&=-a_1(T^1_{11})^2-(\hat{d}_{12}\lambda_2+b_1)T^1_{11}T^1_{21}\\
g^1_1&=-2a_2T^1_{21}T^2_{21}-(\hat{d}_{21}\lambda_1+b_2)(T^1_{11}T^2_{21}+T^2_{11}T^1{21}),&
g^1_2&=-a_2(T^1_{21})^2-(\hat{d}_{21}\lambda_2+b_2)T^1_{11}T^1_{21}
\end{align*}
Studying the cubic truncated dynamical systems, it can be seen that it admits an equilibrium with a Hopf instability, and a necessary condition is $A_1B_1<0$.

Then we can evaluate the necessary condition for the parameter set in Table \ref{tab:param}. The Matlab scripts are available at~\cite{GitHubSL}. We observe that it holds until a certain value of $d_{21}$ corresponding to the value at which the negative region of $L(\lambda_1)$ on the neutral stability curve change position, namely when it appear for $d_{12}>\hat{d}_{12}$. At the same time, when the necessary condition is satisfied, we numerically detected a Hopf bifurcation only for $d_{12}<\hat{d}_{12}$, which produce stable time-periodic patterns. However, when it is not satisfied, the continuation software \texttt{pde2path} detects Hopf bifurcations both for $d_{12}<\hat{d}_{12}$ and $d_{12}>\hat{d}_{12}$, but time-periodic solutions turn out to be unstable. This suggests the presence of a higher co-dimension bifurcation point.

\section{Concluding remarks}\label{sec_concl}
In this work we have analysed the full cross-diffusion SKT model, namely with a cross-diffusion term in both equations, studying the influence of the additional interspecific-competition pressure (cross-diffusion) on the bifurcation structure in the weak-competition regime. This extends the study carried out in \cite{breden2021influence}. In particular, we focused on the type of pitchfork bifurcation on the homogenous branch (related to the stability of the bifurcating branches close to the homogeneous one) and on the presence of Hopf bifurcation points on the bifurcating branch corresponding to the 1-mode.

The model present multistability of solutions, when a particular stable inhomogeneous solution can coexist with the homogeneous one for suitable parameter values. This is an important aspect in ecology, since in this case by perturbing the system it is possible to pass from homogeneous distribution of the species on the habitat to spatial segregation. From the mathematical point of view, this situation is related to the type of pitchfork bifurcations (sub- or super-critical) on the homogenous branch, and its dependence on the cross-diffusion coefficients. In particular, we have here obtained an analytical characterization of the bifurcation point (sub- or super-critical) through weakly nonlinear analysis, deriving the Stuart--Landau equation at the bifurcation point, exploiting the technique presented in \cite{gambino2012turing}. Even though the expression of the coefficient that characterises the pitchfork is nasty and its dependence on the cross-diffusion coefficients is not evident, we can compute it varying the cross-diffusion parameters. As predicted by the linearised analysis, the second cross-diffusion coefficient moves the bifurcation points towards zero, squeezing the bifurcation structure and making it disappear. At the same time, it triggers the appearance of a multi-stability region, where the homogeneous steady state coexists with 1- and/or 2-modes. Note that this effect cannot be captured using linearised analysis only.

From the numerical investigation of the bifurcation structure, we found the appearance of Hopf bifurcation points, suggesting the formation of time-periodic spatial patterns. The influence of the cross-diffusion terms on the Hopf points, as well as on the effective presence, type and stability properties of these time-varying patterns, are biologically relevant. To investigate the possible scenarios, the analytical results obtained in \cite{izuhara2018spatio} have been partially extended to the full cross-diffusion case, and combined with the weakly nonlinear analysis and the numerical continuation. On the one hand the additional cross-diffusion term ``can move'' the system closer to the doubly degenerate point, where a Hopf bifurcation and stable time-period spatial pattern may appear. On the other hand, this doubly degenerate centre seems to change when the additional cross-diffusion coefficient increases, probably because of a higher co-dimension bifurcation. A Hopf bifurcation point is detected, yet the time-periodic solution originated seems to be unstable.

Thanks to its interplay between linearised analysis, weakly nonlinear analysis and numerical continuation, this work constitute a step forward in the analytical understanding of the bifurcation structure of the SKT system, it points out new interesting aspects and opens several different questions that can be addressed in future works.

First of all, a stronger characterisation of the doubly degenerate point at the critical value is needed at this point. This would also allow progress in the understanding of the time-periodic spatial patterns which potentially originate in the vicinity of the doubly degenerate point. While the analytical approach may be feasible, the continuation software \texttt{pde2path} is not immediately suited for the detection of codimension-2 bifurcation points, so this will be a matter of future investigations.

On the other side, it would also be interesting to investigate how far from the doubly degenerate point and why the Hopf bifurcation point and of the time-periodic spatial patterns disappear. From the ecological viewpoint, the influence of the domain size on the type of stable steady and time-periodic patterns is crucial. It is not clear if the domain size has an influence only on the solution profiles (due to different unstable modes), or it can even induce major deformations of the bifurcation structure.
Another important direction is the study of the strong competition case. In \cite{breden2021influence} it has been shown an interesting effect of the additional cross-diffusion term on the bifurcation structure and the presence of Hopf bifurcation points. Better characterization and a deeper investigation would improve the understanding of this different regime.  Note that the derivation of the Stuart--Landau equation and the weakly nonlinear analysis holds and it can predict the type of pitchfork bifurcation on the homogeneous branch also in the strong-competition regimes.

Moreover, an extremely interesting and actual research direction is the extension of cross-diffusion induced instability on networks \cite{EDCKCS}.

Finally, the same study could be carried out for other quasilinear problems involving
cross-diffusion terms. For instance, in the context of predator--prey systems, it is possible to derive by time-scale arguments a different type of cross-diffusion \cite{conforto2018, desvillettes2018}. The linearised analysis suggests that they do not increase the parameter region in which patterns appear, but as in the present work, the global influence cannot be captured only by the linearised analysis. Taken together, these results will better clarify the role of cross-diffusion terms as the key ingredients in pattern formation.
\section*{Acknowledgments}
CS thanks Valeria Giunta and Hirofumi Izuahara for fruitful discussions about the topic of this paper. Partial support by INdAM-GNFM is gratefully acknowledged by CS. 

\bibliographystyle{plain}
\bibliography{bibliography}

\begin{thebibliography}{10}

\bibitem{Ama88}
H.~Amann.
\newblock Dynamic theory of quasilinear parabolic equations. {I}. {A}bstract
  evolution equations.
\newblock {\em Nonlinear Analysis: Theory, Methods \& Applications},
  12(9):895--919, 1988.

\bibitem{Ama90}
H.~Amann.
\newblock Dynamic theory of quasilinear parabolic equations. {II}.
  {R}eaction--diffusion systems.
\newblock {\em Differential and Integral Equations}, 3(1):13--75, 1990.

\bibitem{breden2021computer}
M.~Breden.
\newblock Computer-assisted proofs for some nonlinear diffusion problems.
\newblock {\em Communications in Nonlinear Science and Numerical Simulation},
  109:106292, 2022.

\bibitem{breden2018existence}
M.~Breden and R.~Castelli.
\newblock Existence and instability of steady states for a triangular
  cross-diffusion system: a computer-assisted proof.
\newblock {\em Journal of Differential Equations}, 264(10):6418--6458, 2018.

\bibitem{breden2021influence}
M.~Breden, C.~Kuehn, and C.~Soresina.
\newblock On the influence of cross-diffusion in pattern formation.
\newblock {\em Journal of Computational Dynamics}, 8(2):213--240, 2021.

\bibitem{breden2013global}
M.~Breden, J.-P. Lessard, and M.~Vanicat.
\newblock Global bifurcation diagrams of steady states of systems of {PDE}s via
  rigorous numerics: a 3-component reaction--diffusion system.
\newblock {\em Acta Applicandae Mathematicae}, 128(1):113--152, 2013.

\bibitem{conforto2018}
F.~Conforto, L.~Desvillettes, and C.~Soresina.
\newblock About reaction--diffusion systems involving the {H}olling-type {II}
  and the {B}eddington--{D}e{A}ngelis functional responses for predator--prey
  models.
\newblock {\em Nonlinear Differential Equations and Applications}, 25(3):24,
  2018.

\bibitem{DevLepMouTre15}
L.~Desvillettes, T.~Lepoutre, A.~Moussa, and A.~Trescases.
\newblock On the entropic structure of reaction-cross diffusion systems.
\newblock {\em Communications in Partial Differential Equations},
  40(9):1705--1747, 2015.

\bibitem{desvillettes2018}
L.~Desvillettes and C.~Soresina.
\newblock Non-triangular cross-diffusion systems with predator--prey reaction
  terms.
\newblock {\em Ricerche di Matematica}, pages 1--20, 2018.

\bibitem{desvillettes2015new}
L.~Desvillettes and A.~Trescases.
\newblock New results for triangular reaction cross diffusion system.
\newblock {\em Journal of Mathematical Analysis and Applications},
  430(1):32--59, 2015.

\bibitem{diamond1975assembly}
J.M. Diamond.
\newblock Assembly of species communities.
\newblock In M.L Cody and Diamond J.M., editors, {\em Ecology and Evolution of
  Communities}, pages 342--444. Cambridge, Mass: Harvard Univ Press, 1975.

\bibitem{dohnal2014pde2path}
T.~Dohnal, J.D.M. Rademacher, H.~Uecker, and D.~Wetzel.
\newblock \texttt{pde2path} 2.0: Multi-parameter continuation and periodic
  domains.
\newblock In {\em Proceedings of the 8th European Nonlinear Dynamics
  Conference, ENOC}, volume 2014, 2014.

\bibitem{ehstand2021numerical}
N.~Ehstand, C.~Kuehn, and C.~Soresina.
\newblock Numerical continuation for fractional {PDE}s: sharp teeth and bloated
  snakes.
\newblock {\em Communications in Nonlinear Science and Numerical Simulation},
  98:105762, 2021.

\bibitem{ei1990pattern}
S.-I. Ei and M.~Mimura.
\newblock Pattern formation in heterogeneous reaction--diffusion--advection
  systems with an application to population dynamics.
\newblock {\em SIAM Journal on Mathematical Analysis}, 21(2):346--361, 1990.

\bibitem{GalGarJun03}
G.~Galiano, M.L. Garz{\'o}n, and A.~J{\"u}ngel.
\newblock Semi-discretization in time and numerical convergence of solutions of
  a nonlinear cross-diffusion population model.
\newblock {\em Numerische Mathematik}, 93(4):655--673, 2003.

\bibitem{gambino2016super}
G.~Gambino, M.C. Lombardo, S.~Lupo, and M.~Sammartino.
\newblock Super-critical and sub-critical bifurcations in a reaction-diffusion
  {S}chnakenberg model with linear cross-diffusion.
\newblock {\em Ricerche di Matematica}, 65(2):449--467, 2016.

\bibitem{gambino2012turing}
G.~Gambino, M.C. Lombardo, and M.~Sammartino.
\newblock Turing instability and traveling fronts for a nonlinear
  reaction--diffusion system with cross-diffusion.
\newblock {\em Mathematics and Computers in Simulation}, 82(6):1112--1132,
  2012.

\bibitem{gambino2013pattern}
G.~Gambino, M.C. Lombardo, and M.~Sammartino.
\newblock Pattern formation driven by cross-diffusion in a 2{D} domain.
\newblock {\em Nonlinear Analysis: Real World Applications}, 14(3):1755--1779,
  2013.

\bibitem{Henry}
D.~Henry.
\newblock {\em Geometric Theory of Semilinear Parabolic Equations}.
\newblock Springer, Berlin Heidelberg, Germany, 1981.

\bibitem{iida2006diffusion}
M.~Iida, M.~Mimura, and H.~Ninomiya.
\newblock Diffusion, cross-diffusion and competitive interaction.
\newblock {\em Journal of Mathematical Biology}, 53(4):617--641, 2006.

\bibitem{IidNimYam18}
M.~Iida, H.~Ninomiya, and H.~Yamamoto.
\newblock A review on reaction--diffusion approximation.
\newblock {\em Journal of Elliptic and Parabolic Equations}, 4(2):565--600,
  2018.

\bibitem{izuhara2018spatio}
H.~Izuhara and S.~Kobayashi.
\newblock Spatio-temporal coexistence in the cross-diffusion competition
  system.
\newblock {\em Discrete \& Continuous Dynamical Systems-S}, 14(3):919, 2021.

\bibitem{izuhara2008reaction}
H.~Izuhara and M.~Mimura.
\newblock Reaction-diffusion system approximation to the cross-diffusion
  competition system.
\newblock {\em Hiroshima Mathematical Journal}, 38(2):315--347, 2008.

\bibitem{Jun10}
A.~J{\"u}ngel.
\newblock Diffusive and nondiffusive population models.
\newblock In {\em Mathematical Modeling of Collective Behavior in
  Socio-Economic and Life Sciences}, pages 397--425. Springer, 2010.

\bibitem{Jun16}
A.~J{\"u}ngel.
\newblock {\em Entropy Methods for Diffusive Partial Differential Equations}.
\newblock Springer, 2016.

\bibitem{kan1993stability}
Y.~Kan-On.
\newblock Stability of singularly perturbed solutions to nonlinear diffusion
  systems arising in population dynamics.
\newblock {\em Hiroshima Mathematical Journal}, 23(3):509--536, 1993.

\bibitem{kishimoto1985spatial}
K.~Kishimoto and H.F. Weinberger.
\newblock The spatial homogeneity of stable equilibria of some
  reaction--diffusion systems on convex domains.
\newblock {\em Journal of Differential Equations}, 58:15--21, 1985.

\bibitem{KuehnBook1}
C.~Kuehn.
\newblock {\em PDE Dynamics: An Introduction}.
\newblock SIAM, 2019.

\bibitem{CKetal6}
C.~Kuehn, N.~Berglund, C.~Bick, M.~Engel, T.~Hurth, A.~Iuorio, and C.~Soresina.
\newblock A general view on double limits in differential equations.
\newblock {\em Physica D: Nonlinear Phenomena}, 431:133105, 2022.

\bibitem{EDCKCS}
C.~Kuehn and C.~Soresina.
\newblock Cross-diffusion induced instability on networks.
\newblock {\em in preparation}.

\bibitem{CKCS}
C.~Kuehn and C.~Soresina.
\newblock Numerical continuation for a fast reaction system and its
  cross-diffusion limit.
\newblock {\em SN Partial Differential Equations and Applications}, 1:7, 2020.

\bibitem{levin1974dispersion}
S.A. Levin.
\newblock Dispersion and population interactions.
\newblock {\em The American Naturalist}, 108(960):207--228, 1974.

\bibitem{LouNi96}
Y.~Lou and W.-M. Ni.
\newblock Diffusion, self-diffusion and cross-diffusion.
\newblock {\em Journal of Differential Equations}, 131(1):79--131, 1996.

\bibitem{LouNiYot04}
Y.~Lou, W.-M. Ni, and S.~Yotsutani.
\newblock On a limiting system in the {L}otka--{V}olterra competition with
  cross-diffusion.
\newblock {\em Discrete \& Continuous Dynamical Systems}, 10(1\&2):435--458,
  2004.

\bibitem{lou2015pattern}
Y.~Lou, W.-M. Ni, and S.~Yotsutani.
\newblock Pattern formation in a cross-diffusion system.
\newblock {\em Discrete \& Continuous Dynamical Systems}, 35(4), 2015.

\bibitem{matano1983pattern}
H.~Matano and M.~Mimura.
\newblock Pattern formation in competition-diffusion systems in nonconvex
  domains.
\newblock {\em Publications of the Research Institute for Mathematical
  Sciences}, 19(3):1049--1079, 1983.

\bibitem{mimura1981stationary}
Masayasu Mimura.
\newblock Stationary pattern of some density-dependent diffusion system with
  competitive dynamics.
\newblock {\em Hiroshima Mathematical Journal}, 11(3):621--635, 1981.

\bibitem{mori2018numerical}
T.~Mori, T.~Suzuki, and S.~Yotsutani.
\newblock Numerical approach to existence and stability of stationary solutions
  to a {SKT} cross-diffusion equation.
\newblock {\em Mathematical Models and Methods in Applied Sciences},
  28(11):2191--2210, 2018.

\bibitem{ni2014existence}
W.-M. Ni, Y.~Wu, and Q.~Xu.
\newblock The existence and stability of nontrivial steady states for {SKT}
  competition model with cross diffusion.
\newblock {\em Discrete \& Continuous Dynamical Systems-A}, 34(12):5271--5298,
  2014.

\bibitem{prufert2014oopde}
U.~Pr{\"u}fert.
\newblock \texttt{OOPDE} - an object oriented approach to finite elements in
  {MATLAB}.
\newblock Quickstart Guide, available at \href{http://www. mathe. tu-freiberg.
  de/nmo/mitarbeiter/uwe-pruefert/software}{http://www. mathe. tu-freiberg.
  de/nmo/mitarbeiter/uwe-pruefert/software}, 2014.

\bibitem{shigesada1979spatial}
N.~Shigesada, K.~Kawasaki, and E.~Teramoto.
\newblock Spatial segregation of interacting species.
\newblock {\em Journal of Theoretical Biology}, 79(1):83--99, 1979.

\bibitem{GitHubSKT}
C.~Soresina.
\newblock Supplementary material.
\newblock Matlab scripts for the bifurcation diagrams at
  \url{https://github.com/soresina/fullSKT}, 2021.
\newblock Accessed March 16, 2021.

\bibitem{GitHubSL}
C.~Soresina.
\newblock Supplementary material.
\newblock Matlab scripts for the Stuart--Landau and Hopf coefficients at
  \url{https://github.com/soresina/fullSKT-SL-H}, 2021.
\newblock Accessed April 22, 2021.

\bibitem{uecker2018hopf}
H.~Uecker.
\newblock Hopf bifurcation and time periodic orbits with \texttt{pde2path} --
  algorithms and applications.
\newblock {\em Communications in Computational Physics}, 25:812--852, 2019.

\bibitem{uecker2021continuation}
H.~Uecker.
\newblock Continuation and bifurcation in nonlinear {PDE}s--{A}lgorithms,
  applications, and experiments.
\newblock {\em Jahresbericht der Deutschen Mathematiker-Vereinigung}, pages
  1--38, 2021.

\bibitem{uecker2021numerical}
H.~Uecker.
\newblock Numerical continuation and bifurcation in nonlinear {PDE}s, 2021.

\bibitem{uecker2014pde2path}
H.~Uecker, D.~Wetzel, and J.D.M. Rademacher.
\newblock \texttt{pde2path} - {A} {M}atlab package for continuation and
  bifurcation in 2{D} elliptic systems.
\newblock {\em Numerical Mathematics: Theory, Methods and Applications},
  7(1):58--106, 2014.

\bibitem{wilson1975}
E.O. Wilson.
\newblock {\em Sociobiology: The New Synthesis}.
\newblock Cambridge: Harvard, 1975.

\bibitem{wollkind1994weakly}
D.J. Wollkind, V.S. Manoranjan, and L.~Zhang.
\newblock Weakly nonlinear stability analyses of prototype reaction–diffusion
  model equations.
\newblock {\em SIAM Reviews}, 36(2):176--214, 1994.

\end{thebibliography}

\end{document}